\documentclass[11pt]{article}
\usepackage[utf8]{inputenc}
\usepackage[a4paper, total={6.25in, 9.25in}]{geometry}
\usepackage{amsmath}
\usepackage{amsthm}
\usepackage{bm}

\usepackage{amssymb}

\usepackage{mathtools}
\usepackage{float}
\usepackage{algorithm}
\usepackage{algpseudocode}

\usepackage{subcaption}
\usepackage{graphicx}
\usepackage{color}

\usepackage[square,comma,numbers,sort&compress]{natbib}
\usepackage[linkcolor=black,colorlinks=true,citecolor=black,filecolor=black,urlcolor=black,pdftitle={A Kernel Framework for Linear Operator Equations},pdfauthor={Damir Filipovic, Valentin Tissot-Daguette}]{hyperref}
\usepackage[capitalise,nameinlink]{cleveref}

\def\be{\begin{equation}}
\def\ee{\end{equation}}

\def\bea{\begin{align}}
\def\eea{\end{align}}

\def\bea*{\begin{align*}}
\def\eea*{\end{align*}}
\theoremstyle{plain}
\newtheorem{theorem}{Theorem}[section]

\newtheorem{lemma}[theorem]{Lemma}
\newtheorem{corollary}[theorem]{Corollary}
\newtheorem{proposition}[theorem]{Proposition}
\newtheorem{asm}[theorem]{Assumption}

\theoremstyle{definition}
\newtheorem{definition}[theorem]{Definition}
\newtheorem{example}[theorem]{Example}
\newtheorem{remark}[theorem]{Remark}

\newcommand{\E}{\mathbb{E}}

\newcommand{\N}{\mathbb{N}}

\newcommand{\Pb}{\mathbb{P}}
\newcommand{\Q}{\mathbb{Q}}
\newcommand{\R}{\mathbb{R}}

\newcommand{\calA}{\mathcal{A}}
\newcommand{\calB}{\mathcal{B}}

\newcommand{\calD}{\mathcal{D}}

\newcommand{\calF}{\mathcal{F}}
\newcommand{\calG}{\mathcal{G}}
\newcommand{\calH}{\mathcal{H}}
\newcommand{\calI}{\mathcal{I}}

\newcommand{\calN}{\mathcal{N}}

\newcommand{\calS}{\mathcal{S}}

\newcommand{\calV}{\mathcal{V}}
\newcommand{\calW}{\mathcal{W}}
\newcommand{\calX}{\mathcal{X}}
\newcommand{\calY}{\mathcal{Y}}

\crefname{asm}{Assumption}{Assumptions}
\Crefname{asm}{Assumption}{Assumptions}
\crefname{figure}{Figure}{Figures}

\usepackage{authblk}
\title{A Kernel Framework for Linear Operator Equations}

\author{Damir Filipovic\thanks{
   EPFL, EXTRA 218, 1015 Lausanne, Switzerland. Email:   damir.filipovic@epfl.ch } \hspace{1em}
Valentin Tissot-Daguette\thanks{
 Quantitative Research, Office of the CTO, Bloomberg, 731 Lexington Avenue, New York, NY 10022, USA. Email: vtissotdague@bloomberg.net}
}

\date{\today}

\begin{document}
\maketitle

\vspace{-4mm}

\begin{abstract}
Solving linear operator equations with nonlocal terms or heterogeneous constraints often requires specific numerical schemes. We develop a unified approach based on embedding such equations into a reproducing kernel Hilbert space (RKHS). Differential, boundary, and integral operators are represented by their Riesz representers in the RKHS. The equations are solved in the regularized least-squares sense at finitely many collocation points. The solution is a finite linear combination of operator-induced representers and is obtained from a single linear system. Its sensitivities, such as financial Greeks, follow by direct differentiation. To reduce the computational cost for large collocation sets, we introduce a Nystr\"om approximation based on subsampling the operator-induced representers. It yields the best approximation of this solution within the reduced trial space, in a regularized empirical norm. 
Numerical examples cover ordinary, partial, and partial integro-differential equations, including multi-asset and jump-diffusion option pricing.
\end{abstract}

{\small
\noindent\textbf{Keywords.} Reproducing kernel Hilbert spaces, linear  operator equations,  least-squares kernel collocation, Nystr\"om approximation, partial integro-differential equations.}
\vspace{3mm}

{\small
\noindent\textbf{MSC Codes.} 
65J10, 
65N35, 
46E22, 
91G60. 
}

\section{Introduction}\label{sec:intro}

Let $\calS$ be a set and $\calF(\calS)$ the space of all real-valued
functions on $\calS$. For $m=1,\ldots,M$, let
\(\calA_m:\calD(\calA_m)\subseteq\calF(\calS)\to\calF(\calX_m)\)
be a linear operator, for some set $\calX_m$, and let
$f_m\in\calF(\calX_m)$. We study systems of linear operator equations
\begin{equation}\label{eq:generalProblem}
    \calA_m u=f_m
    \quad\text{on }\calX_m,
    \qquad m=1,\ldots,M.
\end{equation}
This formulation covers ordinary and partial differential equations (ODEs and PDEs), integral equations, and partial
integro-differential equations (PIDEs). For a PDE, for instance, one operator
describes the dynamics in the interior of the domain. Another operator encodes a
boundary condition through a restriction operator. Integral and other nonlocal
terms are included in the same way.

Numerical methods for \eqref{eq:generalProblem} are typically tailored to the
structure of the operators. Classical finite-difference,
finite-element, and spectral methods work well when this structure
can be exploited. They become harder to use in high dimensions, for
nonlocal operators, or when several heterogeneous constraints
must be handled at the same time.

Kernel methods provide a mesh-free alternative. In a reproducing kernel
Hilbert space (RKHS), continuous linear functionals admit Riesz representers.
Differential, boundary, and integral operators can thus be represented
in a common framework. This observation underlies classical
generalized interpolation and kernel collocation methods, and more
recent Gaussian-process approaches to differential equations.

\paragraph{Contributions.}
We develop such a framework by embedding \eqref{eq:generalProblem} into an RKHS
$\calH\subseteq\calF(\calS)$.
First, differential, boundary, and integral operators are represented by their
Riesz representers in $\calH$. We give sufficient conditions for this
representability for differential operators and for nonlocal jump operators of
L\'evy type. This leads to a regularized population
least-squares problem and a normal equation in the
RKHS.

Second, we evaluate the equations at finitely many collocation points. The
solution of the resulting problem, which we call the kernel solution, is a finite linear combination of operator-induced representers. This follows from the
classical representer theorem for bounded linear functionals. Its coefficients solve a linear system with a block structure, in which each pair $(\calA_\ell,\calA_m)$ induces its own cross-kernel.
Sensitivities of the solution are obtained by differentiating its kernel representation, provided the kernel is sufficiently smooth. In financial
applications this yields, for example, Delta and Gamma without solving additional
equations. The kernel derivatives are evaluated by
automatic differentiation.

Third, for large collocation sets we introduce a Nystr\"om approximation. It
restricts the trial space to the span of a subset of the operator-induced
representers. The reduced solution solves the same empirical problem over this smaller space. For any choice of the subspace, it is the best approximation of the kernel solution within this subspace, in a regularized empirical norm.
We also bound the expected gap in the objective under uniform random subsampling.

We illustrate the framework on ODEs, parabolic PDEs, and pricing PIDEs with jump
terms. The examples include multidimensional option-pricing equations. They show
how nonlocal terms are included through the operator-induced
representers. In higher dimensions we use low-discrepancy Sobol sequences for the
collocation points.

\paragraph{Related literature.}
Our approach builds on the theory of reproducing kernel Hilbert
spaces of Aronszajn~\cite{aronszajn1950theory}. See also
Paulsen and Raghupathi~\cite{paulsen2016}. Reproducing properties of derivative and integral operators
in an RKHS are studied in Zhou~\cite{Zhou}, Steinwart and Christmann~\cite[Corollary~4.36]{SteinwartChristmann2008},
and El-Boukkouri et al.~\cite{ElBoukkouri2025}. The representer theorem for regularized problems
involving general bounded linear functionals goes back to
Kimeldorf and Wahba~\cite{KimeldorfWahba1971}. See also Wahba~\cite[Chapter~1]{Wahba1990}.
Generalized interpolation by linear functionals and its
connection with mesh-free approximation and PDE collocation are classical.
See, for example, Wendland~\cite[Chapter~16]{Wendland2005},
Schaback and Wendland~\cite{SchabackWendland2006}, and Fasshauer~\cite{Fasshauer2007}.
Early radial-basis-function collocation methods include the multiquadric scheme
of Kansa \cite{KANSA1990127,KANSA1990147}. Using the operator-induced
representers as trial functions is known as symmetric collocation
\cite{Fasshauer1997,FrankeSchaback1998}. Our kernel solution is a regularized
least-squares variant of it. Least-squares kernel collocation and
applications to PIDEs and stochastic PDEs are discussed, for example, in
\cite{Cheung,ChanHubbert2010,Cialenco01122012}. A kernel-based approach to solving path-dependent PDEs was also recently proposed by Pannier and Salvi~\cite{PannierSalvi}.

A related viewpoint interprets kernel collocation through
Gaussian processes and optimal recovery. Graepel \cite{Graepel2003} formulates the
solution of noisy linear operator equations as Gaussian-process regression with
operator-transformed kernels. Pf\"ortner et al.  \cite{Pfoertner2022} show that
Gaussian-process regression with observations through bounded linear operators
generalizes classical solvers such as collocation and Galerkin methods.
Chen et al.~\cite{ChenHosseiniOwhadiStuart2021} develop a kernel and Gaussian-process framework for
solving nonlinear PDEs and inverse problems.
Batlle et al.~\cite{BatlleChenHosseiniOwhadiStuart2025} establish quantitative convergence
rates for nonlinear and parametric PDEs. The last two works mainly use hard
collocation constraints. Our formulation starts from a regularized
least-squares problem for several, possibly heterogeneous, linear
operators, including nonlocal ones. We also consider its population counterpart.

The action of differential operators on Gaussian-process covariance kernels has
also been used in more data-driven settings. For example, Raissi et al.~\cite{RAISSI2017}
infer parameters of prescribed linear differential operators from scarce and
noisy observations. In our setting, the operators
and forcing terms are known, and the goal is to solve a given system of
operator equations.

Our goal also differs from that of the literature on \emph{operator learning}. Methods such as DeepONet
\cite{LuJinPang2021}, Fourier Neural Operators
\cite{LiKovachkiAzizzadenesheli2021}, and neural operators more generally
\cite{KovachkiLiLiu2023} learn maps between function spaces from
input--output pairs. They typically aim at approximating the solution
operator for a family of PDEs. Kernel-based approaches to operator
learning have also been developed. See, for example,
Batlle et al.~\cite{BatlleDarcyHosseiniOwhadi2024}. In the present work, the
operators $\calA_m$ and forcing terms $f_m$ are prescribed, and the goal is to
solve a single system \eqref{eq:generalProblem}. Our method is thus a numerical
solver and not a surrogate model learned across a family of problems.

Finally, Nystr\"om approximations are widely used to reduce the computational
cost of large kernel systems
\cite{WilliamsSeeger,RudiCamorianoRosasco2015,AlaouiMahoney2015,filipovic2025kernel}.
The standard learning-theory results quantify the quality of a randomly selected
subspace with high probability. We use the Nystr\"om approximation as a numerical tool for fixed collocation points. For any selected subset, the reduced solution satisfies a deterministic best-approximation property. Random subsampling is one way of selecting the~subset.

\paragraph{Outline.}
\Cref{sec:prelim} reviews the RKHS preliminaries and the representation of linear
operators. \Cref{sec:linearOperatorEqs} develops the population and empirical problems,
the representation of the kernel solution, the estimation of sensitivities, and the Nystr\"om
approximation. \Cref{sec:numerics} introduces a kernel adapted to parabolic scaling and presents numerical examples for the Bachelier and multi-asset Black--Scholes PDEs, as well as pricing PIDEs with nonlocal jump terms.  \Cref{sec:conclusion} concludes.
The proofs are collected in Appendix~\ref{app:proofs}, while additional numerical examples, including ordinary differential equations (ODEs) and the forward heat equation, are presented in Appendix~\ref{app:suppNumerics}. 

\section{Preliminaries}\label{sec:prelim}

Let $\calS$ be a set, and denote by $\calF(\calS)$ the space of all real-valued functions on $\calS$. Let $\calH\subseteq \calF(\calS)$ be an RKHS with reproducing kernel $k:\calS\times\calS \to \R$. Denote the inner product and norm in  $\calH$ by $( \cdot,\cdot )_{\calH}$ and $\|h\|_{\calH} = \sqrt{( h,h )_{\calH}}$, respectively. The kernel $k$ has the reproducing property
\begin{equation}\label{eq:reproduce}
    h(x) = ( k_x,h )_{\calH},\quad \forall  h \in \calH,  \quad x\in \calS,
\end{equation}
where the \emph{kernel section} $k_x \coloneqq k(\cdot,x) \in\calH$ is the Riesz \emph{representer} of the  evaluation at~$x$. Throughout, $|\cdot|$ denotes absolute value and $\|\cdot\|_2$ the Euclidean norm on $\R^d$.

\subsection{Reproducing Linear Operators}\label{sec:reproducingOperator}
Consider a linear operator  $\calA:\calD(\calA)\subseteq\calF(\calS)\to \calF(\calX)$ for some set $\calX$. Assume that the RKHS $\calH$ is contained in the domain $\calD(\calA)$ of $\calA$:
\begin{equation}\label{HinDA}
    \calH\subseteq \calD(\calA).
\end{equation}
In particular,  $k_y\in\calD(\calA)$ and  $\calA k_y\in\calF(\calX)$ for all $y\in \calS$. We can thus define the \emph{operator-induced kernel}
\begin{equation}\label{kcalAdef}
 k^\calA:\calX\times\calS\to\R,\quad    k^\calA(x,y) \coloneqq \calA k_y(x),
\end{equation}
and we write $k^\calA_{x} \coloneqq k^\calA(x,\cdot) \in\calF(\calS)$ for $x\in\calX$.
Our framework requires that operators like $\calA$ can be represented within $\calH$. The following property makes this precise.
\begin{definition}\label{def:ARegular}
We say that $\calH$ is \emph{$\calA$-regular} if \eqref{HinDA} holds and, for each $x\in\calX$, the linear functional $h\mapsto\calA h(x)$ is continuous on $\calH$.
\end{definition}

The Riesz representer of this functional is necessarily given by the operator-induced kernel.

\begin{lemma}\label{lem:ARegular}
Assume that \eqref{HinDA} holds. Then $\calH$ is $\calA$-regular if and only if, for each $x\in\calX$, the function $k^\calA_x$ lies in $\calH$ and the reproducing property
\begin{equation}\label{reprAh}
    \calA h(x) = ( k^\calA_{x},h )_{\calH}
\end{equation}
holds for all $h\in\calH$.
\end{lemma}

We call $k^\calA_x$ the \emph{representer of $\calA$ at $x$}, and refer to the functions $k^\calA_x$, $x\in\calX$, as \emph{operator-induced representers}.
Regularity is stable under linear combinations with variable coefficients. Indeed, let $\calH$ be $\calA$- and $\calB$-regular for two linear operators $\calA,\calB$ with values in $\calF(\calX)$, and let $a,b\in\calF(\calX)$. Then $\calH$ is regular for the operator $a\calA+b\calB$, defined by $(a\calA+b\calB)h(x)\coloneqq a(x)\calA h(x)+b(x)\calB h(x)$, with representers
\begin{equation}\label{eq:lincomb}
    k^{a\calA+b\calB}_x = a(x)\,k^\calA_x+b(x)\,k^\calB_x,\quad x\in\calX.
\end{equation}

We start with an elementary example.
\begin{example}\label{extraceop}
Let $\calX\subseteq\calS$ and $\calA h=h|_\calX$ be the restriction (or trace) operator with domain $\calD(\calA)=\calF(\calS)$. By \eqref{eq:reproduce}, $\calA h(x)=h(x)=(k_x,h)_\calH$ for all $h\in\calH$ and $x\in\calX$. Hence $\calH$ is $\calA$-regular,  with representer $k^\calA_x=k_x$.
\end{example}
The following result covers differential operators. For an open set $\calX\subseteq\R^d$ with $\calX\subseteq\calS$, and $1\le i,j\le d$, we denote by $\partial_i h(x)$ the partial derivative of $h\in\calF(\calS)$ at $x\in\calX$ in the $i$-th coordinate direction. The domain $\calD(\partial_i)$ consists of all $h$ for which $\partial_ih(x)$ exists for every $x\in\calX$, and we write $\partial_{ij}h\coloneqq\partial_i(\partial_jh)$, with the corresponding domain. For the kernel, $\partial_{x_i}$ and $\partial_{y_i}$ denote the partial derivatives with respect to the $i$-th coordinate of the first and second argument, respectively.

\begin{proposition}\label{propdiff}
 Let $\calX\subseteq\R^d$ be an open set with $\calX\subseteq\calS$, and let $1\le i, j\le d$.
 \begin{enumerate}\renewcommand{\theenumi}{(\roman{enumi})}\renewcommand{\labelenumi}{\theenumi}
 \item\label{propdiff1} Assume that the partial derivatives $\partial_{y_i}k$ and $\partial_{x_i}\partial_{y_i}k$ exist and are continuous on $\calX\times\calX$. Then $\calH$ is $\partial_i$-regular, with
 \begin{equation}\label{eqkpartdiff1}
     k_x^{\partial_i}(y)=\partial_{x_i}k(x,y),  \quad 
     \|k_x^{\partial_i}\|_\calH^2 =\partial_{x_i y_i}k(x,y)\big|_{y=x},  \quad  x\in\calX,\ y\in\calS. 
 \end{equation}
 \item\label{propdiff2} Assume that the hypothesis of \ref{propdiff1} holds for the index $j$, and that the function $\widetilde k\coloneqq\partial_{x_j}\partial_{y_j}k$ on $\calX\times\calX$ satisfies the hypothesis of \ref{propdiff1} for the index $i$. Then $\calH$ is $\partial_{ij}$-regular, with
 \begin{equation}\label{eqkpartdiff}
     k_x^{\partial_{ij}}(y)=\partial_{x_ix_j}k(x,y),  
     \|k_x^{\partial_{ij}}\|_\calH^2=\partial_{x_i y_i}\widetilde k(x,y)\big|_{y=x}, \quad  x\in\calX,\ y\in\calS. 
 \end{equation}
 \end{enumerate}
\end{proposition}

The assumptions of \cref{propdiff} hold for all $i,j$ if the restriction of $k$ to $\calX\times\calX$ is four times continuously differentiable. They are local. Neither boundedness of $k$ and its derivatives nor compactness of $\calX$ is required. For compact domains, the result is due to Zhou~\cite[Theorem~1]{Zhou}. For open domains, see also Steinwart and Christmann~\cite[Corollary~4.36]{SteinwartChristmann2008}.

\begin{figure}
\caption{Gaussian kernel section $k_{x}$ and representers $k^{\partial_1}_{x}, k^{\partial_{11}}_{x}$ as functions of $y$, for $\gamma=2$ and $x\in\{-1.5,0,1.5\}$ (dots).}
        \centering
        \includegraphics[width=0.98\linewidth]{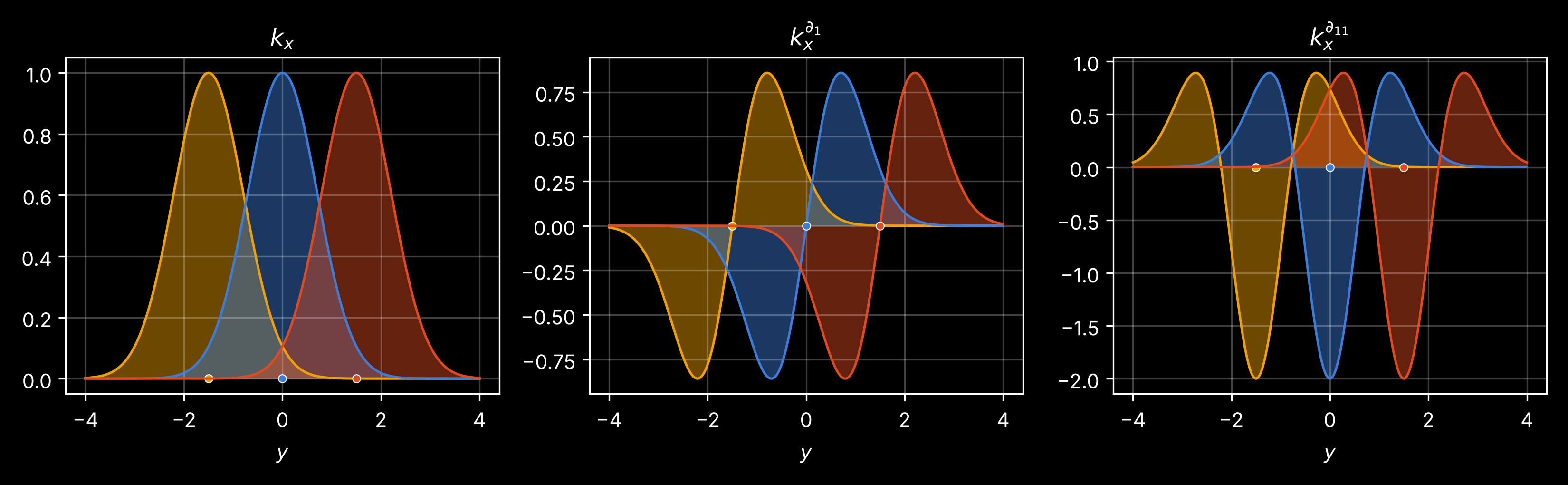}
    \label{fig:GaussianKernel}
\end{figure}

\begin{example} \label{ex:GaussianKernel}
      Let $\calS = \R$,  and consider the Gaussian kernel $k(x,y) = e^{-\frac{\gamma}{2}|x-y|^2}$, $\gamma>0$. A direct calculation shows that    \begin{align}\label{eq:DerivativesGaussian}
        k_x^{\partial_1}(y) &= \partial_xk(x,y)= - \gamma (x-y) k_x(y), \\[1em]
        k_x^{\partial_{11}}(y) &= \partial_{xx}k(x,y)= (\gamma^2 (x -y)^2  - \gamma  )k_x(y).
    \end{align}
  Figure~\ref{fig:GaussianKernel} displays the kernel section $k_x$ and the representers $k_x^{\partial_{1}}$, $k_x^{\partial_{11}}$.
\end{example}

\begin{example}\label{ex:integralEquations}
Let $\calS$ be compact and consider the Fredholm equation
\[
\int_{\calS}u(y)g(x,y)\,dy=f(x),\qquad x\in\calS,
\]
with continuous $g:\calS\times\calS\to\R$. If $k$ is continuous on $\calS\times\calS$, then the associated RKHS $\calH$ is regular for the operator
$\calA u(x)=\int_{\calS}u(y)g(x,y)\,dy$, with representer
\[
k_x^\calA=\int_{\calS}k_y\,g(x,y)\,dy.
\]
Indeed, this follows directly by interchanging the RKHS inner product with the Bochner integral.
\end{example}

The same argument applies to the nonlocal operators that arise in jump models, where the integral is taken with respect to a possibly infinite measure.

\begin{proposition}\label{prop:jump}
Let $\calS\subseteq\R^d$ be a Borel set, $\calX\subseteq\calS$, and let $\mu$ be a Borel measure on $\R^d\setminus\{0\}$ such that
\begin{equation}\label{eq:finiteVariation}
    \int \min(1,\|z\|_2)\,\mu(dz)<\infty,
\end{equation}
and $x+z\in\calS$ for all $x\in\calX$ and $\mu$-almost all $z$. Assume that $k$ is continuous and that there exist constants $c_0,c_1\ge 0$ such that
\begin{equation}\label{eq:jumpAssumption}
    k(y,y)\le c_0^2\quad\text{and}\quad \|k_{x+z}-k_x\|_\calH\le c_1\|z\|_2,
\end{equation}
for all $y\in\calS$, $x\in\calX$, and $\mu$-almost all $z$. Then $\calH$ is $\calA$-regular for the operator
\begin{equation}\label{eq:jumpOperator}
    \calA u(x)\coloneqq\int\big(u(x+z)-u(x)\big)\,\mu(dz),\quad x\in\calX,
\end{equation}
whose domain $\calD(\calA)$ consists of all $u\in\calF(\calS)$ for which the integrand is $\mu$-integrable for every $x\in\calX$. Its representers are given by the Bochner integrals
\begin{equation}\label{eq:jumpRepresenter}
\begin{aligned}
    k^\calA_x=\int\big(k_{x+z}-k_x\big)\,\mu(dz),\qquad 
    \|k^\calA_x\|_\calH
    \le\int\min(2c_0,c_1\|z\|_2)\,\mu(dz),
    \qquad x\in\calX.
\end{aligned}
\end{equation}
\end{proposition}

Condition \eqref{eq:finiteVariation} covers finite measures, as in the case of compound Poisson jumps. It also covers L\'evy measures of infinite activity and finite variation, such as that of the Variance Gamma model in \cref{sec:VG}. For the Gaussian kernel $k(x,y)=e^{-\frac\gamma2\|x-y\|_2^2}$ on $\R^d$, condition \eqref{eq:jumpAssumption} holds with $c_0=1$ and $c_1=\sqrt\gamma$. Note that $\calS$ must contain $x+z$ for all $x\in\calX$ and $\mu$-almost all $z$. If $\mu$ has unbounded support, $\calS$ is therefore unbounded in the jump directions, even when the set $\calX$ is bounded.

\subsection{\texorpdfstring{$L^2(\Pb)$}{L2(P)} Embeddings}

Let $\calS$ be a set and $(\calX,\Pb)$ be a probability space. Let $\calA:\calD(\calA)\subseteq\calF(\calS)\to\calF(\calX)$ be a linear operator. Let $\calH$ be an RKHS on $\calS$ with reproducing kernel $k$. Following the usual convention, we identify $\Pb$-square integrable functions with their $\Pb$-a.s.\ equivalence classes in $L^2(\Pb)$.

\begin{lemma}\label{lem:L2}
Assume that $\calH$ is separable and $\calA$-regular, that the function $x\mapsto k^\calA(x,y)$ is measurable on $\calX$ for every $y\in\calS$, and that
\begin{equation}\label{assL2}
    \int_{\calX} \| k^\calA_x\|_\calH^2\,\Pb(dx)<\infty.
\end{equation}
Then the following hold:
\begin{enumerate}
    \item The map $x\mapsto k^\calA_x:\calX\to\calH$ is measurable, and for every $g\in L^2(\Pb)$ the Bochner integral $\int_{\calX} k^\calA_x\,g(x)\,\Pb(dx)$ is a well-defined element of $\calH$.
    \item $\calA h$ is a $\Pb$-square integrable function on $\calX$ for all $h\in\calH$.
    \item $\calA :\calH\to L^2(\Pb)$ is Hilbert--Schmidt, with 
    $ \|\calA\|_\text{\textnormal{HS}}^2 = \int_{\calX} \| k^\calA_x\|_\calH^2\,\Pb(dx).$

    \item The adjoint $\calA^\ast:L^2(\Pb)\to\calH$ of $\calA$ satisfies
    \begin{equation}\label{eq:adjoint}
        \calA^\ast g =\int_{\calX} k^\calA_x\, g(x)\,\Pb(dx).
    \end{equation}
\end{enumerate}
\end{lemma}

\section{Linear Operator Equations}\label{sec:linearOperatorEqs}

We now return to the system of linear operator equations \eqref{eq:generalProblem}. Let $(\calX_m,\Pb_m)$ be probability spaces, $m=1,\dots,M$. We seek a solution of \eqref{eq:generalProblem} within a separable RKHS $\calH$ on $\calS$ with reproducing kernel $k$, and  assume that $\calH$ is  $\calA_m$-regular for all $m=1,\dots,M$. \cref{ex:zhou} below provides a prototype of such a system. In view of \cref{lem:ARegular}, equation \eqref{eq:generalProblem} can  be rewritten as
\begin{equation}
    ( k^{\calA_m}_{x},h )_{\calH} =f_m(x) , \quad \forall x \in \calX_m, \;  m = 1,\ldots, M,
\end{equation}
where $h\in\calH$ is the unknown. If $\calA_m$ is the restriction operator (see \cref{extraceop}), then the $m$-th equation reduces to $h(x)=( k_{x},h )_{\calH} = f_m(x)$ for all $x \in \calX_m$.
\begin{example}\label{ex:zhou}
  Given a domain $\Omega \subseteq \R^d$, $d\in \N$, and functions $f\in\calF(\Omega)$ and $g\in\calF(\partial\Omega)$, consider the Dirichlet problem
  \begin{equation}
      \begin{cases}
 \calG u = f, \quad & \textnormal{on \; } \Omega, \\
     u  = g, \quad & \textnormal{on \;} \partial\Omega,
 \end{cases}
 \label{eq:Dirichlet}
  \end{equation}
with the differential operator
    $ \calG =  \sum_{i,j = 1}^d a_{ij}\partial_{ij} + \sum_{i=1}^d b_i \partial_i  + c$, for given coefficient functions  $a_{ij}, b_i,c\in  \calF(\Omega)$.
  Then problem \eqref{eq:Dirichlet}  can be rewritten as  \eqref{eq:generalProblem} by setting $M=2$ and
\begin{align*}
\calS=\overline\Omega,\quad
\calX_1=\Omega,\quad
\calX_2=\partial\Omega,\quad 
f_1=f,\quad f_2=g,\quad
\calA_1u=\calG u,\quad
\calA_2u=u|_{\partial\Omega}.
\end{align*}
If the restriction of $k$ to $\Omega\times\Omega$ is four times continuously differentiable, it follows from \cref{propdiff} and \eqref{eq:lincomb} that $\calH$ is $\calA_1$-regular, with representers
$k^{\calA_1}_x=\sum_{i,j}a_{ij}(x)k^{\partial_{ij}}_x+\sum_i b_i(x)k^{\partial_i}_x+c(x)k_x$. For the restriction operator $\calA_2$, the representers are given by Example~\ref{extraceop}.
\end{example}

\subsection{Population Problem}\label{sec:population}
Given objective weights $w_m >0$ and a Tikhonov regularization parameter $\lambda >0$, we consider the regularized population problem
\begin{equation}\label{eq:population}
\begin{aligned}
&\underset{h\in\calH}{\operatorname{minimize}}\quad \Phi_\lambda(h),\\
&\Phi_\lambda(h)\coloneqq
\sum_{m=1}^M w_m\|\calA_mh-f_m\|_{L^2(\Pb_m)}^2
+\lambda\|h\|_\calH^2.
\end{aligned}
\end{equation}
\begin{asm}\label{asm:L2}
For all $m=1,\ldots,M$, the function $x\mapsto k^{\calA_m}(x,y)$ is measurable on $\calX_m$ for every $y\in\calS$,
\begin{equation}
   \int_{\calX_m} \| k^{\calA_m}_x\|_{\calH}^2 \, \Pb_m(dx) <\infty,
\end{equation}
and $f_m\in L^2(\Pb_m)$.
\end{asm}
By \cref{lem:L2}, \cref{asm:L2} implies that $\calA_m:\calH\to L^2(\Pb_m)$ is Hilbert--Schmidt, with adjoint $\calA_m^*$ given by \eqref{eq:adjoint}. The integrability condition holds in particular if the representers are bounded,
\begin{equation}\label{eq:kappa}
    \kappa_m\coloneqq\sup_{x\in\calX_m}\|k^{\calA_m}_x\|_\calH<\infty.
\end{equation}
Problem \eqref{eq:population} is related to \eqref{eq:generalProblem} as follows. If \eqref{eq:generalProblem} admits a solution $u\in\calH$, then $\Phi_\lambda(u)=\lambda\|u\|_\calH^2$. Note that \eqref{eq:population} only involves \eqref{eq:generalProblem} $\Pb_m$-almost everywhere.

We introduce the compact notation
\begin{gather*}
\bm{\calA}h = (\calA_1h,\ldots, \calA_Mh)^{\top}, \qquad  \bm{f}= (f_1,\ldots, f_M)^{\top}, \\
\bm{\calA}^*\bm f = \sum_{m=1}^M\calA_m^*f_m, \qquad \bm{\calW}\bm{f} = ( w_1 f_1  ,\ldots,  w_M f_M )^\top,
\end{gather*}
so that $\bm\calA:\calH\to L^2(\Pb_1)\times\cdots\times L^2(\Pb_M)$ has adjoint $\bm\calA^*$, and we denote by $\bm{I}$ the identity on $\calH$.
\begin{proposition}\label{prop:population}
   Let \cref{asm:L2} hold. Then, for all $\lambda>0$, the population problem \eqref{eq:population} admits a unique solution $h_\lambda\in\calH$. It is given by
   \begin{equation}\label{eq:hLambda}
        h_\lambda =   (\bm{\calA}^* \bm{\calW} \bm{\calA} +\lambda \bm{I})^{-1}  \bm{\calA}^*  \bm{\calW} \bm{f},
    \end{equation}
    or, equivalently,
\begin{equation}\label{eq:hLambda2}
        h_\lambda = \bm{\calA}^*(\bm{\calA}\bm{\calA}^* + \lambda\bm{\calW}^{-1})^{-1} \bm{f}.
    \end{equation}
\end{proposition}
The second expression \eqref{eq:hLambda2} is the one used for computations, see \cref{sec:finiteSample}.
For all $\ell, m \le M$,  we have
\begin{align*}
\calA_{\ell}\calA_m^*f(x)
=\int_{\calX_m}k^{\ell,m}(x,y)\,f(y)\,\Pb_m(dy), \qquad x\in\calX_\ell,\; f\in L^2(\Pb_m).
\end{align*}
with the \emph{cross-kernel} $k^{\ell,m}$ on $\calX_\ell \times  \calX_m$ induced by the pair $(\calA_\ell, \calA_{m})$,
\begin{equation}\label{eq:mlkernel}
\begin{aligned}
k^{\ell,m}(x,y)
\coloneqq (k^{\calA_\ell}_x,k^{\calA_m}_y)_\calH,\qquad  x\in\calX_\ell,\; y\in\calX_m.
\end{aligned}
\end{equation}

\begin{remark}[The limit $\lambda\to0$]\label{rem:lbda0} Standard Tikhonov theory implies that, whenever the unregularized least-squares problem admits a minimizer, \(h_\lambda\) converges as \(\lambda\to0\) to its minimum-norm solution; see \cite[Chapters~2--5]{EnglHankeNeubauer1996}.
\end{remark}

\subsection{Empirical Problem}
\label{sec:finiteSample}
Let $\bm{x}_m = (x_{m,1},\ldots,x_{m,n_m}) \in \calX_m^{n_m}$ be given collocation points, $m = 1,\ldots, M$. They may be drawn independently from $\Pb_m$, or be generated by a low-discrepancy sequence or a regular grid. In what follows they are fixed. Consider the sampling operators and forcing vectors
\begin{align*}
    &A_m: \calH \to \R^{n_m}, \quad A_m h
    = (\calA_m h(x_{m,1}),\ldots,\calA_m h(x_{m,n_m}))^{\top}, \\[1em]
    &\bm{y}_m
    = (f_m(x_{m,1}),\ldots,f_m(x_{m,n_m}))^{\top}\in \R^{n_m}.
\end{align*}
Replacing each $\Pb_m$ in \eqref{eq:population} by the empirical measure $\frac1{n_m}\sum_{n=1}^{n_m}\delta_{x_{m,n}}$ of the collocation points yields the \emph{empirical problem}
\begin{equation}\label{eq:SampleObjective}
\begin{aligned}
&\underset{h\in\calH}{\operatorname{minimize}}\quad \widehat\Phi_\lambda(h),\\
&\widehat\Phi_\lambda(h)\coloneqq
\sum_{m=1}^M\frac{w_m}{n_m}\|A_mh-\bm y_m\|_2^2
+\lambda\|h\|_\calH^2.
\end{aligned}
\end{equation}
We write $N\coloneqq\sum_{m=1}^Mn_m$ for the total number of collocation points, and
\begin{align*}
\bm{A}&=(A_1,\ldots,A_M)^\top:\calH\to\R^N,\\
\bm y&=(\bm y_1,\ldots,\bm y_M)^\top,\\
\bm{G}&=\operatorname{diag}\Big(
\frac{n_1}{w_1}I_{n_1},\ldots,
\frac{n_M}{w_M}I_{n_M}\Big).
\end{align*}
The adjoint operator $A_m^*:\R^{n_m}\to \calH$ with respect to the Euclidean inner product
maps a vector of coefficients $\beta \in \R^{n_m}$ to a linear combination of operator-induced representers,
   \begin{equation}\label{eq:adjointEmp}
A_m^* \beta =
\sum_{n=1}^{n_m}\beta_n \ k^{\calA_m}_{x_{m,n}}  \in \calH.
\end{equation}
The action of $A_\ell$ on these adjoint operators yields the $(\ell,m)$-block
of the $N\times N$ matrix $\bm{A}\bm{A}^*$. Its entries are pairwise evaluations of the
cross-kernel $k^{\ell,m}$ in \eqref{eq:mlkernel}:
\begin{equation}\label{eq:blocks}
(A_\ell A_m^*)_{n,n'}
=
k^{\ell,m}(x_{\ell,n},x_{m,n'})
=
\big(
\calA_\ell\calA_m k
\big)
(x_{\ell,n},x_{m,n'}), \quad  n\le n_\ell, \, n'\le n_m,
\end{equation}
where $\calA_\ell$ and $\calA_m$ act on the first and
second arguments of $k$, respectively.
Arguing as for \cref{prop:population}, and using \eqref{eq:adjointEmp}, we obtain the following representation.

\begin{proposition} \label{thm:representer}%
    For all $\lambda>0$, the empirical problem \eqref{eq:SampleObjective} admits a unique solution $\widehat h_\lambda\in\calH$. It is given by
    \begin{equation}\label{eq:hLambda3}
        \widehat h_\lambda =\bm{A}^*\widehat{\bm\beta}_\lambda,\quad 
        \widehat{\bm{\beta}}_\lambda
        =(\widehat\beta_1,\ldots,\widehat\beta_M)
        =(\bm{A}\bm{A}^*+\lambda\bm{G})^{-1}\bm{y} 
        \in\R^{n_1}\times\cdots\times\R^{n_M}.
    \end{equation}
    that is,
    \begin{equation}\label{eq:representer}
    \widehat h_\lambda(x)
    =\sum_{m=1}^M\sum_{n=1}^{n_m}
    \widehat\beta_{m,n}\,k^{\calA_m}_{x_{m,n}}(x),
    \qquad x\in\calS.
    \end{equation}
\end{proposition}

We refer to $\widehat h_\lambda$ as the \emph{kernel solution}. \cref{thm:representer} is an instance of the classical representer theorem for regularized problems with bounded linear functionals \cite{KimeldorfWahba1971,Wahba1990}. See also Wendland~\cite[Chapter~16]{Wendland2005} for the interpolation limit $\lambda\to0$, which corresponds to symmetric collocation.

The kernel solution also admits the usual Gaussian-process posterior-mean interpretation with observation variances \(\lambda n_m/w_m\); see \cite{Graepel2003,Pfoertner2022}.

\begin{remark}[Empirical versus population problem]
We do not analyze the convergence
of $\widehat h_\lambda$ to $h_\lambda$, or to a solution of
\eqref{eq:generalProblem}, as the number of collocation points grows. For hard-constrained
kernel collocation, quantitative rates in terms of the fill distance are given
in Batlle et al.~\cite{BatlleChenHosseiniOwhadiStuart2025}. For regularized least squares with
randomly drawn points, see Blanchard and M{\"u}cke~\cite{BlanchardMuecke2018} in the context of inverse problems and Rudi et al.~\cite{RudiCamorianoRosasco2015} for regression.
\end{remark}

\subsection{Algorithm and Implementation Details}\label{rem:numericalDetails}

We call the resulting method the \emph{kernel solver}. Its training and prediction phases are summarized in \cref{alg:KernelSolver}.
When $\calA_m$ is a differential operator, we compute the operator-induced kernels by automatic differentiation with \texttt{JAX}~\cite{jax2018github}. For the second-order PDE examples below, the blocks \eqref{eq:blocks} require fourth-order mixed derivatives, making automatic differentiation particularly convenient. Integral operators are evaluated in closed form when available, as in \cref{sec:jumpdiffusion}, and otherwise by numerical quadrature, as in \cref{sec:VG}. 

\begin{algorithm}[H]
\caption{Kernel Solver for Linear Operator Equations}\label{alg:KernelSolver}

\vspace{1mm}
\textbf{Input:} Problem data $(\calA_m,f_m,\calX_m)_{m=1}^M$;
collocation sizes $(n_m)_{m=1}^M$; kernel $k:\calS\times\calS\to\R$;
regularization parameter $\lambda>0$; weights $(w_m)_{m=1}^M$.\\[0.1em]
\textbf{Output:} Coefficient vectors $(\widehat\beta_m)_{m=1}^M$ and kernel solution $\widehat h_\lambda\in\calH$.
\\[-0.8em]
\noindent\rule{\textwidth}{0.5pt}
\textbf{\underline{TRAINING}}
\begin{itemize}
\setlength\itemsep{0.5em}

\item[I.] \textbf{Generate} collocation points $x_{m,n}\in\calX_m$,
$n=1,\ldots,n_m$, $m=1,\ldots,M$.

\item[II.] \textbf{Compute}
\begin{enumerate}
    \item \textbf{Operator blocks:} $(\bm{A}\bm{A}^*)_{\ell,m}
    =\big(k^{\ell,m}(x_{\ell,n},x_{m,n'})\big)
    \in\R^{n_\ell\times n_m},$

    \item \textbf{Forcing vectors:} $    \bm{y_m}
    =\big(f_m(x_{m,1}),\ldots,f_m(x_{m,n_m})\big)^\top
    \in\R^{n_m}.$
\end{enumerate}

\item[III.] \textbf{Assemble} $\bm{A}\bm{A}^*+\lambda\bm{G}$, where
\(
\bm{G}=\operatorname{diag}\Big(
\frac{n_1}{w_1}I_{n_1},\ldots,
\frac{n_M}{w_M}I_{n_M}\Big).
\)

\item[IV.] \textbf{Return}
\((\widehat\beta_1,\ldots,\widehat\beta_M)
=(\bm{A}\bm{A}^*+\lambda\bm{G})^{-1}\bm y.\)

\end{itemize}

\vspace{2mm}

\textbf{\underline{PREDICTION}}
\begin{itemize}
    \item[I.] \textbf{Given} a query point $x\in\calS$.
    \item[II.] \textbf{Evaluate}
    $k^{\calA_m}_{x_{m,n}}(x)=\calA_m k(x_{m,n},x)$
    for $n=1,\ldots,n_m$ and $m=1,\ldots,M$.
    \item[III.] \textbf{Return}
    \[
    \widehat h_\lambda(x)
    =\sum_{m=1}^M\sum_{n=1}^{n_m}
    \widehat\beta_{m,n}\,k^{\calA_m}_{x_{m,n}}(x).
    \]
\end{itemize}
\end{algorithm}

\subsection{Estimating Sensitivities}\label{sec:Greeks}

Many applications, such as the hedging of financial derivatives, require the solution of \eqref{eq:generalProblem} and also its sensitivities with respect to the state variables. By the representation \eqref{eq:representer}, these sensitivities are obtained by applying the corresponding operator to the operator-induced representers. See also \cite{LiuLi} for a related approach in regression problems.

\begin{corollary}\label{prop:sensitivities}
Let $\calB:\calD(\calB)\subseteq\calF(\calS)\to\calF(\calY)$ be a linear operator, for some set $\calY$, such that $\calH$ is $\calB$-regular. Then, for all $x\in\calY$, the kernel solution $\widehat{h}_\lambda$ in \eqref{eq:representer} satisfies
\begin{align}\label{eq:sensB}
&\calB\widehat h_\lambda(x)
=\sum_{m=1}^M\sum_{n=1}^{n_m}
\widehat\beta_{m,n}\,\calB k^{\calA_m}_{x_{m,n}}(x), \\[0.5em]
\label{eq:sensBound}
    &|\calB h(x)-\calB h'(x)|\le\|k^\calB_x\|_\calH\,\|h-h'\|_\calH,\quad h,h'\in\calH.
\end{align}
In particular, if $\calS \subseteq \R^{d}$ and the assumptions of \cref{propdiff} hold on an open set $\calY\subseteq\calS$, then, for $x\in\calY$,
     \begin{align}
\partial_i\widehat h_\lambda(x)
=\sum_{m=1}^M\sum_{n=1}^{n_m}
\widehat\beta_{m,n}\,\partial_i k^{\calA_m}_{x_{m,n}}(x),
\quad 
\partial_{ij}\widehat h_\lambda(x)
=\sum_{m=1}^M\sum_{n=1}^{n_m}
\widehat\beta_{m,n}\,\partial_{ij} k^{\calA_m}_{x_{m,n}}(x).
\label{eq:derivatives}
     \end{align}
\end{corollary}

Indeed, \eqref{eq:sensB} follows from \eqref{eq:representer} and the linearity of $\calB$, and \eqref{eq:sensBound} from \eqref{reprAh} and the Cauchy--Schwarz inequality. By \eqref{eq:sensBound}, convergence in $\calH$ implies pointwise convergence of all sensitivities for which $\calH$ is regular.
The kernel derivatives are computed by automatic differentiation, see Section~\ref{rem:numericalDetails}.

\begin{example} \label{ex:Greeks}
    Let $u(t,x) = \E^{\Q}[g(X_T) \ | \  X_t = x]$ be the price of a European-style derivative with   payoff $g:\R\to \R$ and  maturity $T>0$. We consider the Bachelier model with zero interest rates and volatility $\sigma>0$, that is, $X_t = X_0 + \sigma W_t$ under the risk-neutral measure $\Q$, for some $\Q$-Brownian motion $W$. Then $u$ solves the parabolic PDE
\begin{equation}\label{eq:pricingPDE}
\begin{aligned}
\partial_tu+\frac12\sigma^2\partial_{xx}u&=0,
&&\textnormal{on }\calX_1\coloneqq(0,T)\times\R,\\
u&=g,
&&\textnormal{on }\calX_2\coloneqq\{T\}\times\R.
\end{aligned}
\end{equation}
  As in Example~\ref{ex:zhou}, we set $f_1 = 0$, $f_2 = g$, $\calA_1 = \partial_t + \frac{1}{2}\sigma^2 \partial_{xx}$, and let $\calA_2$ be the restriction~operator.

The sensitivities of $u$ with respect to the state variables are called the \emph{Greeks}. They are used in the risk management of financial derivatives. The \emph{Delta} $\Delta \coloneqq \partial_x u$ is the change in the option price per unit change in the underlying. It is used for dynamic hedging. If $\widehat{h}_{\lambda}$ denotes the kernel solution of \eqref{eq:pricingPDE} and $(t_{1,n},x_{1,n})_{n=1}^{n_1}$, $(T,x_{2,n})_{n=1}^{n_2}$ are collocation points in $\calX_1$ and $\calX_2$, respectively, then \cref{prop:sensitivities}  provides the following  Delta estimate
   \begin{equation}
   \begin{aligned}
\widehat\Delta_\lambda(t,x)
&=\sum_{n=1}^{n_1}\widehat\beta_{1,n}\,
\partial_x k^{\calA_1}\big((t_{1,n},x_{1,n}),(t,x)\big) +
\sum_{n=1}^{n_2}\widehat\beta_{2,n}\,
\partial_x k\big((T,x_{2,n}),(t,x)\big).
   \end{aligned}
   \end{equation}
   Numerical results are presented in \cref{sec:Bachelier}.
\end{example}


\subsection{Nystr\"om Approximation}\label{sec:Nystrom}

Computing the kernel solution $\widehat h_\lambda$ in \eqref{eq:hLambda3}
requires solving a linear system of size $N$.
When $N$ is large, the computational cost can be reduced by restricting the empirical problem \eqref{eq:SampleObjective} to a lower-dimensional subspace of $\calH$. We use the span of a subset of the operator-induced representers in \eqref{eq:representer}.

For each $m=1,\ldots,M$, let
\(\mathcal{J}_m\subseteq\{1,\ldots,n_m\}\), 
\(|\mathcal{J}_m|=r_m\ge1,
\) 
be an index set, and set $r\coloneqq\sum_{m=1}^M r_m$.
We define the Nystr\"om subspace
\begin{equation}\label{eq:NystromSpace}
\calH_r
\coloneqq
\operatorname{span}
\left\{
k^{\calA_m}_{x_{m,n}}
:
n\in\mathcal{J}_m,\ m=1,\ldots,M
\right\}
\subseteq\calH,
\end{equation}
and the \emph{reduced solution}
\begin{equation}\label{eq:NystromEstimator}
\widehat h_{\lambda,r}
\coloneqq
\underset{h\in\calH_r}{\operatorname{argmin}}\,
\widehat\Phi_\lambda(h).
\end{equation}
Hence $\widehat h_{\lambda,r}$ solves the same empirical problem as the kernel solution $\widehat h_\lambda$, but over a smaller trial space spanned by $r$ selected operator-induced representers.

Relabel the selected representers as
$\varphi_1,\ldots,\varphi_r$ and write
\[
\widehat h_{\lambda,r}
=
\sum_{j=1}^r\widehat\alpha_j\varphi_j.
\]
Define
\[
(K_r)_{j,j'}
\coloneqq
(\varphi_j,\varphi_{j'})_{\calH},
\qquad
(Z_m)_{n,j}
\coloneqq
\calA_m\varphi_j(x_{m,n}),
\]
so that $K_r\in\R^{r\times r}$ and
$Z_m\in\R^{n_m\times r}$. The coefficient vector
$\widehat\alpha\in\R^r$ satisfies
\begin{equation}\label{eq:NystromSystem}
\left(
\sum_{m=1}^M
\frac{w_m}{n_m}Z_m^\top Z_m
+\lambda K_r
\right)\widehat\alpha
=
\sum_{m=1}^M
\frac{w_m}{n_m}Z_m^\top\bm y_m.
\end{equation}
If the selected representers are linearly dependent, the coefficient vector
$\widehat\alpha$ need not be unique. The corresponding element
$\widehat h_{\lambda,r}\in\calH_r$ is nevertheless unique, since $\lambda>0$. A coefficient vector can be obtained using the Moore--Penrose pseudoinverse. The reduced problem requires the $N\times r$ cross-kernel evaluations in $Z_1,\ldots,Z_M$ instead of the full $N\times N$ matrix $\bm{A}\bm{A}^*$. Forming the system \eqref{eq:NystromSystem} requires of the order of $Nr^2$ operations. It remains to solve an $r$-dimensional~system.

\paragraph{Best approximation.}
We first quantify the error induced by the restriction to a subspace. Define the empirical normal operator
\[
\widehat T
\coloneqq
\sum_{m=1}^M\frac{w_m}{n_m}A_m^*A_m,
\]
so the kernel solution satisfies the normal equation $(\widehat T+\lambda \bm I)\widehat h_\lambda=\sum_{m=1}^M\frac{w_m}{n_m}A_m^*\bm y_m$. The operator $\widehat T$ is bounded, self-adjoint, and positive semidefinite on $\calH$. Hence
\[(h,h')_{\widehat T,\lambda}\coloneqq((\widehat T+\lambda \bm I)h,h')_\calH\]
defines an inner product on $\calH$, with the regularized empirical norm
\begin{equation}\label{eq:empiricalNorm}
\|h\|_{\widehat T,\lambda}^2
=
(h,\widehat T h)_{\calH}
+\lambda\|h\|_{\calH}^2
=
\sum_{m=1}^M
\frac{w_m}{n_m}\|A_mh\|_2^2
+\lambda\|h\|_{\calH}^2,
\end{equation}
which is equivalent to the norm of $\calH$,
\begin{equation}\label{eq:normEquivalence}
\lambda\|h\|_{\calH}^2\le\|h\|_{\widehat T,\lambda}^2\le\big(\|\widehat T\|+\lambda\big)\|h\|_{\calH}^2 ,
\end{equation}
where $\|\widehat T\|$ denotes the operator norm. We observe that  \eqref{eq:empiricalNorm} is the energy norm of the quadratic functional $\widehat\Phi_\lambda$. Such norms are standard in the analysis of Galerkin methods. There, the best-approximation property below is known as C\'ea's lemma, see Ciarlet~\cite{Ciarlet1978}. Closely related norms, induced by regularized covariance operators, are used in the analysis of regularized least squares and of Nystr\"om approximations in learning theory, see, e.g., Rudi et al.~\cite{RudiCamorianoRosasco2015}.

\begin{proposition}\label{prop:Nystrom}
Let $\calV\subseteq\calH$ be a closed subspace and 
$\widehat h_{\lambda,\calV}\coloneqq\operatorname{argmin}_{h\in\calV}\widehat\Phi_\lambda(h)$. Then $\widehat h_{\lambda,\calV}$ is the orthogonal projection of the kernel solution $\widehat h_\lambda$ onto $\calV$ with respect to the inner product $(\cdot,\cdot)_{\widehat T,\lambda}$. In particular,
\begin{align}
\|\widehat h_{\lambda,\calV}-\widehat h_\lambda\|_{\widehat T,\lambda}
&=
\min_{h\in\calV}
\|h-\widehat h_\lambda\|_{\widehat T,\lambda},
\label{eq:NystromBestApprox}
\\[0.5em]
\widehat\Phi_\lambda(\widehat h_{\lambda,\calV})
-
\widehat\Phi_\lambda(\widehat h_\lambda)
&=
\|\widehat h_{\lambda,\calV}-\widehat h_\lambda\|_{\widehat T,\lambda}^2,
\label{eq:NystromObjectiveGap}
\end{align}
and
\begin{align}
\|\widehat h_{\lambda,\calV}-\widehat h_\lambda\|_{\calH}
&\le
\frac{1}{\sqrt{\lambda}}
\|\widehat h_{\lambda,\calV}-\widehat h_\lambda\|_{\widehat T,\lambda},
\label{eq:NystromHBound}\\[0.5em]
\left(
\sum_{m=1}^M
\frac{w_m}{n_m}
\|A_m(\widehat h_{\lambda,\calV}-\widehat h_\lambda)\|_2^2
\right)^{1/2}
&\le
\|\widehat h_{\lambda,\calV}-\widehat h_\lambda\|_{\widehat T,\lambda}.
\label{eq:NystromResidualBound}
\end{align}
\end{proposition}

Proposition~\ref{prop:Nystrom} applies in particular to the Nystr\"om subspace $\calV=\calH_r$, for any choice of the index sets $\mathcal J_m$. It shows that the only error of the Nystr\"om approximation comes from the restriction of the trial space. The reduced solution $\widehat h_{\lambda,r}$ is the best approximation of the kernel solution
$\widehat h_\lambda$ in the regularized empirical norm. We refer to the left-hand side of \eqref{eq:NystromObjectiveGap}, with $\calV=\calH_r$, as the \emph{objective gap} of the reduced solution. Combining \eqref{eq:NystromHBound} with \eqref{eq:sensBound} gives a pointwise bound on the error of the reduced solution and of its sensitivities. For every operator $\calB$ as in \cref{prop:sensitivities},
\begin{equation}
|\calB\widehat h_{\lambda,r}(x)-\calB\widehat h_\lambda(x)|\le
\frac{\|k^\calB_x\|_\calH}{\sqrt\lambda}\,
\|\widehat h_{\lambda,r}-\widehat h_\lambda\|_{\widehat T,\lambda}.
\end{equation}

\paragraph{Random subsampling.}
We now assume that the index sets are chosen at random and bound the expected objective gap. The collocation points remain fixed. The only randomness is the choice of the index sets $(\mathcal J_m)_{m=1}^M$. We denote by $\E_{\mathcal J}$ the expectation with respect to it.

\begin{proposition}\label{prop:NystromRate}
Assume $n_m\ge2$ and let
$\mathcal{J}_m$ be drawn uniformly at random among all subsets of $\{1,\ldots,n_m\}$ of size $r_m$, with $1\le r_m\le n_m$, independently across $m$.
Write
\[
\varphi_{m,n}\coloneqq k^{\calA_m}_{x_{m,n}},
\qquad
\widehat h_{\lambda,m}
\coloneqq
\sum_{n=1}^{n_m}\widehat\beta_{m,n}\varphi_{m,n},
\qquad\text{so that}\qquad
\widehat h_\lambda
=
\sum_{m=1}^M\widehat h_{\lambda,m},
\]
and define
\[
V_{\lambda,m}
\coloneqq
\frac1{n_m}\sum_{n=1}^{n_m}
\left\|
n_m\widehat\beta_{m,n}\varphi_{m,n}
-\widehat h_{\lambda,m}
\right\|_{\widehat T,\lambda}^2.
\]
Then
\begin{equation}\label{eq:NystromRandomRate}
\E_{\mathcal J}\!\left[
\widehat\Phi_\lambda(\widehat h_{\lambda,r})
-\widehat\Phi_\lambda(\widehat h_\lambda)
\right]
\le
\sum_{m=1}^M
\frac{n_m-r_m}{r_m(n_m-1)}
V_{\lambda,m}.
\end{equation}
Moreover,
\begin{equation}\label{eq:VmExplicit}
V_{\lambda,m}
=
n_m\sum_{n=1}^{n_m}
\widehat\beta_{m,n}^2
\|\varphi_{m,n}\|_{\widehat T,\lambda}^2
-
\|\widehat h_{\lambda,m}\|_{\widehat T,\lambda}^2.
\end{equation}
\end{proposition}

The quantities $V_{\lambda,m}$ can be bounded in terms of the residuals of the kernel solution.

\begin{corollary}\label{cor:NystromExplicit}
Under the assumptions of Proposition~\ref{prop:NystromRate}, define
\[
\kappa_{\lambda,m,N}^2
\coloneqq
\max_{1\le n\le n_m}
\|\varphi_{m,n}\|_{\widehat T,\lambda}^2,
\qquad
R_m^2
\coloneqq
\frac{w_m}{n_m}\,
\|\bm y_m-A_m\widehat h_\lambda\|_2^2,
\qquad
R_N^2
\coloneqq
\sum_{m=1}^M R_m^2.
\]
Then
\begin{align}
\E_{\mathcal J}\!\left[
\widehat\Phi_\lambda(\widehat h_{\lambda,r})
-\widehat\Phi_\lambda(\widehat h_\lambda)
\right]
&\le
\frac1{\lambda^2}
\sum_{m=1}^M
\frac{n_m-r_m}{r_m(n_m-1)}
\,w_m\kappa_{\lambda,m,N}^2R_m^2
\nonumber\\[0.5em]
&\le
\frac{R_N^2}{\lambda^2}
\max_{1\le m\le M}
\left\{
w_m\kappa_{\lambda,m,N}^2
\frac{n_m-r_m}{r_m(n_m-1)}
\right\}.
\label{eq:NystromResidualRate}
\end{align}
Moreover, the following hold:
\begin{enumerate}\renewcommand{\theenumi}{(\roman{enumi})}\renewcommand{\labelenumi}{\theenumi}
\item\label{cor:NystromExplicit1} The residual satisfies
$R_N^2
\le
\widehat\Phi_\lambda(\widehat h_\lambda)
\le
\widehat\Phi_\lambda(0)
=
\sum_{m=1}^M
\frac{w_m}{n_m}\|\bm y_m\|_2^2$. If \eqref{eq:generalProblem} admits a solution $u\in\calH$, then also $R_N^2\le\lambda\|u\|_\calH^2$.
\item\label{cor:NystromExplicit2} If the representers are bounded, $\kappa_m<\infty$ in \eqref{eq:kappa} for all $m$, then
\[\kappa_{\lambda,m,N}^2\le\kappa_m^2\Big(\sum_{\ell=1}^Mw_\ell\kappa_\ell^2+\lambda\Big).\]
\end{enumerate}
\end{corollary}

For equal numbers of collocation points and of selected representers across equations, we obtain the following bound.

\begin{corollary}\label{cor:NystromBalancedRate}
Under the assumptions of Proposition~\ref{prop:NystromRate}, suppose that
$n_m=n$ and $r_m=r_0$ for all $m\le M$, so that $r=Mr_0$. Then
\begin{align}
\E_{\mathcal J}\!\left[
\widehat\Phi_\lambda(\widehat h_{\lambda,r})
-\widehat\Phi_\lambda(\widehat h_\lambda)
\right]
&\le
\frac{n-r_0}{r_0(n-1)}\,
\frac{R_N^2}{\lambda^2}
\max_{1\le m\le M}
\left\{
w_m\kappa_{\lambda,m,N}^2
\right\}\nonumber\\[0.5em]
&\le
\frac{M}{r}\,
\frac{R_N^2}{\lambda^2}
\max_{1\le m\le M}
\left\{
w_m\kappa_{\lambda,m,N}^2
\right\}.
\label{eq:NystromBalancedRate}
\end{align}
If, in addition, $\kappa_m<\infty$ and $f_m$ is bounded for all $m$, then the right-hand side of \eqref{eq:NystromBalancedRate} is bounded by $C/(\lambda^2r)$, with the constant
\[
C= M\,
\Big(\sum_{m=1}^M w_m\sup_{\calX_m}|f_m|^2\Big)
\max_{1\le m\le M}\{w_m\kappa_m^2\}\times
\Big(\sum_{\ell=1}^M w_\ell\kappa_\ell^2+\lambda\Big).
\]
which depends neither on $N$ nor on the collocation points. If, moreover, \eqref{eq:generalProblem} admits a solution $u\in\calH$, then the factor $R_N^2/\lambda^2$ in \eqref{eq:NystromBalancedRate} can be replaced by $\|u\|_\calH^2/\lambda$.
\end{corollary}

The bound \eqref{eq:NystromBalancedRate} is of order $r^{-1}$ in the number of selected representers. Its constant is uniform in the number of collocation points when the representers and the forcing terms are bounded. The bound is conservative, for two reasons. First, it is obtained by comparing the reduced solution with a reweighted subsample of the kernel solution. This is in general a weaker competitor than the projection in \cref{prop:Nystrom}. Second, its constant scales like $R_N^2/\lambda^2$, or like $\|u\|_\calH^2/\lambda$ if \eqref{eq:generalProblem} is solvable in $\calH$. The quantities $V_{\lambda,m}$ in \eqref{eq:VmExplicit} can be computed from the kernel solution and give a sharper estimate.

\section{Numerical Examples}\label{sec:numerics}

Our examples concern backward parabolic partial differential equations (PDEs) and partial integro-differential equations (PIDEs). 
The collocation points are restricted to a bounded set in the spatial variables whenever the underlying domain is unbounded; this set is a box except in \cref{sec:maxcall}. Equivalently, the sets $\calX_m$ below are replaced by their intersection with this bounded set. The kernels used below are smooth, so that \cref{propdiff} applies without further conditions. On a bounded box, the coefficients of the equations and the relevant kernel derivatives are bounded. The representers of the jump operators in \cref{sec:jumpdiffusion,sec:VG} are bounded by \eqref{eq:jumpRepresenter}. Hence the representers are bounded as in \eqref{eq:kappa}, and \cref{asm:L2} and the bounds of \cref{sec:Nystrom} apply. 
Further implementation details are given in Section~\ref{rem:numericalDetails}. 

 For diffusion-driven problems, it is useful to adapt the kernel to the parabolic scaling in time to maturity. Given
$(t,x),(s,y)\in\calS\coloneqq[0,T]\times\R^d$, consider
\begin{align}
k\big((t,x),(s,y)\big)
&=
\exp\left\{
-\frac{\gamma_0}{2}|t-s|^2
-\frac{\gamma_1}{2}\|\psi(t,x)-\psi(s,y)\|_2^2
\right\},
\label{eq:parabolicKernel}\\[0.5em]
\psi(t,x)
&=\frac{x}{\sqrt{T-t+\varepsilon}},
\qquad \varepsilon>0.
\label{eq:transformBwd}
\end{align}
The transformation $\psi$ reflects the parabolic scaling in time to maturity: at fixed $t$, the spatial length scale is proportional to
$\sqrt{T-t+\varepsilon}$. The parameter $\varepsilon>0$ prevents a singularity at maturity and ensures that the kernel remains smooth.

\subsection{Bachelier Pricing PDE}
\label{sec:Bachelier}
 Consider an at-the-money call option with payoff $g(x)=(x-x_0)^+$ in the Bachelier model of \cref{ex:Greeks}, with $X_0=x_0$. The price function $u(t,x) = \E^{\Q}[(X_T - x_0)^+ \  | \ X_t = x ]$ solves the backward heat equation \eqref{eq:pricingPDE}.
We use the parabolic kernel \eqref{eq:parabolicKernel} with spatial transformation $\psi$ given in \eqref{eq:transformBwd}. We set $\gamma_0=8$ and $\gamma_1 = \frac{1}{\sigma^2}$, so that the spatial length scale of the kernel matches the diffusion scale $\sigma\sqrt{T}$ of the Bachelier model. We take $T=1$, $\sigma=20$, $x_0=100$, and the spatial box $[x_0-6\sigma\sqrt T,\,x_0+6\sigma\sqrt T]=[-20,220]$. The collocation points are drawn uniformly in $(0,T)\times[-20,220]$ and in $\{T\}\times[-20,220]$. We use $\lambda=10^{-7}$, $w_1=1$, and $w_2=10$ to put greater emphasis on the terminal condition.

\begin{figure}
\caption{Price and sensitivities (Delta, Gamma) of a call option in the Bachelier model. Solid lines: kernel solution. Dashed lines: exact solution.}
    \centering
      \begin{subfigure}[t]{0.325\linewidth}
        \centering
        \caption{Price $(u)$}
        \includegraphics[width = 0.99\linewidth, height = 1.6in]{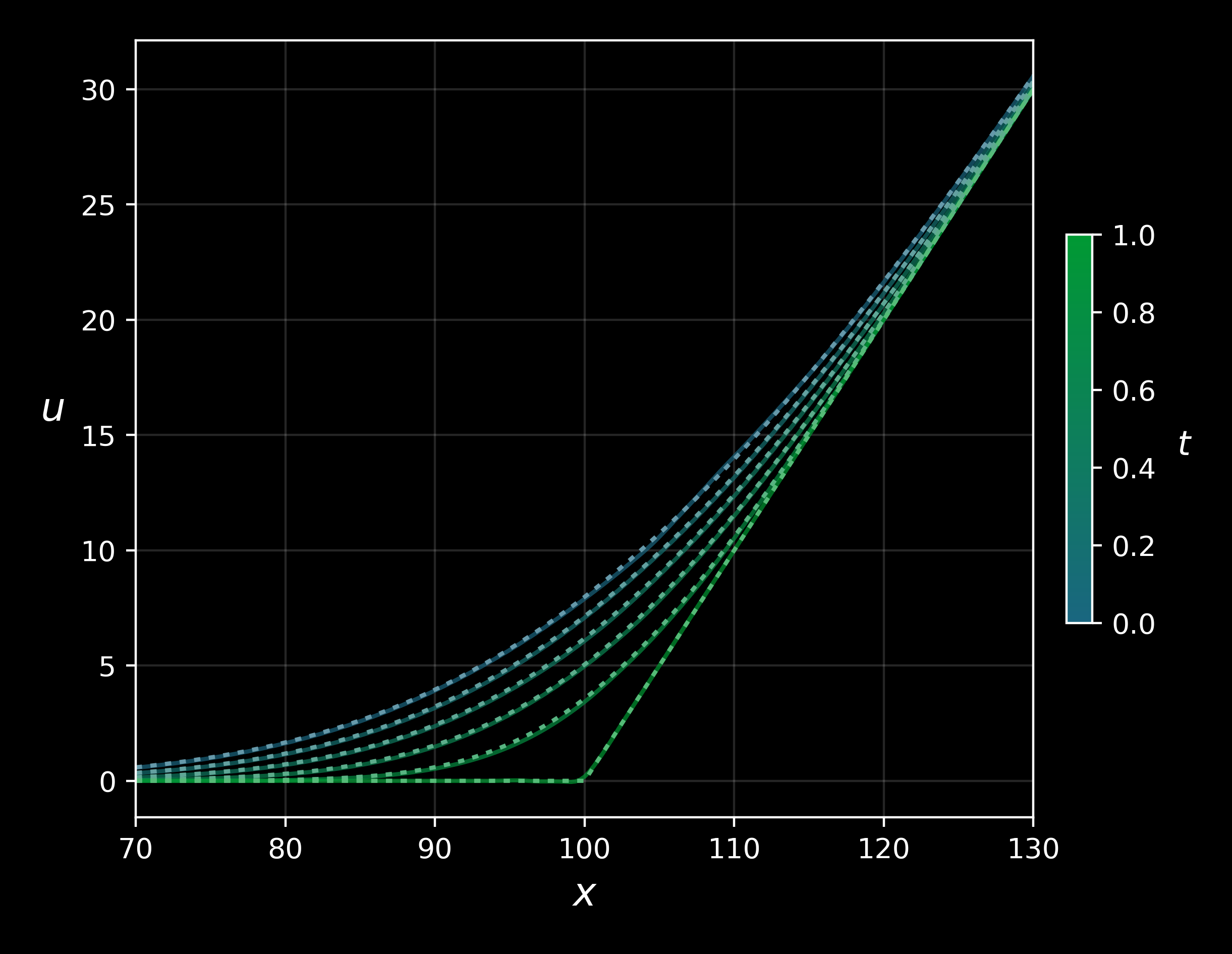}
        \label{fig:BachelierSolution}
    \end{subfigure}
       \hfill
    \begin{subfigure}[t]{0.325\linewidth}
        \centering
       \caption{Delta $(\Delta)$}
        \includegraphics[width = 0.99\linewidth, height = 1.6in]{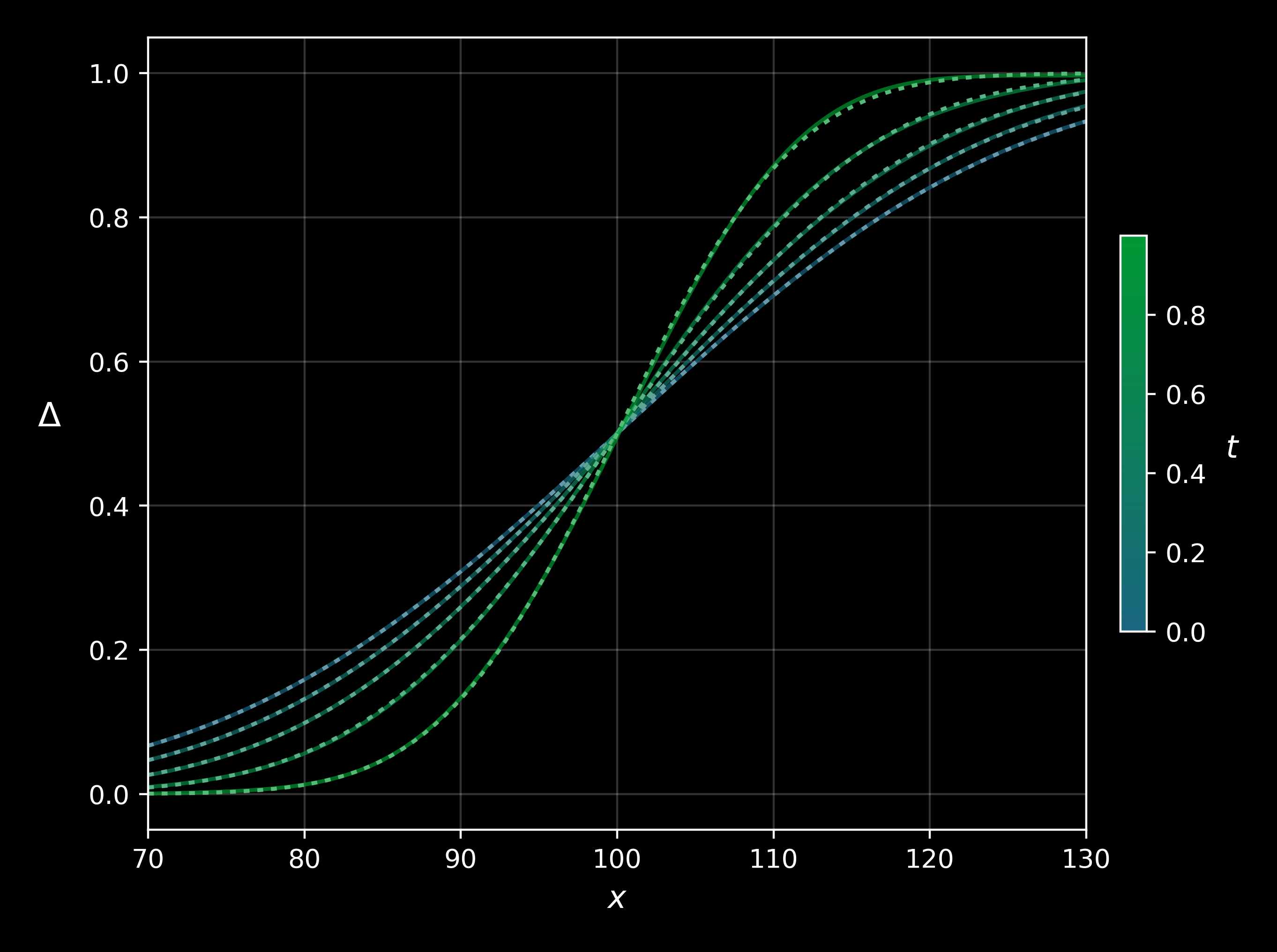}
        \label{fig:BachelierDelta}
    \end{subfigure}
        \begin{subfigure}[t]{0.325\linewidth}
        \centering
        \caption{Gamma $(\Gamma)$}
        \includegraphics[width = 0.99\linewidth, height = 1.6in]{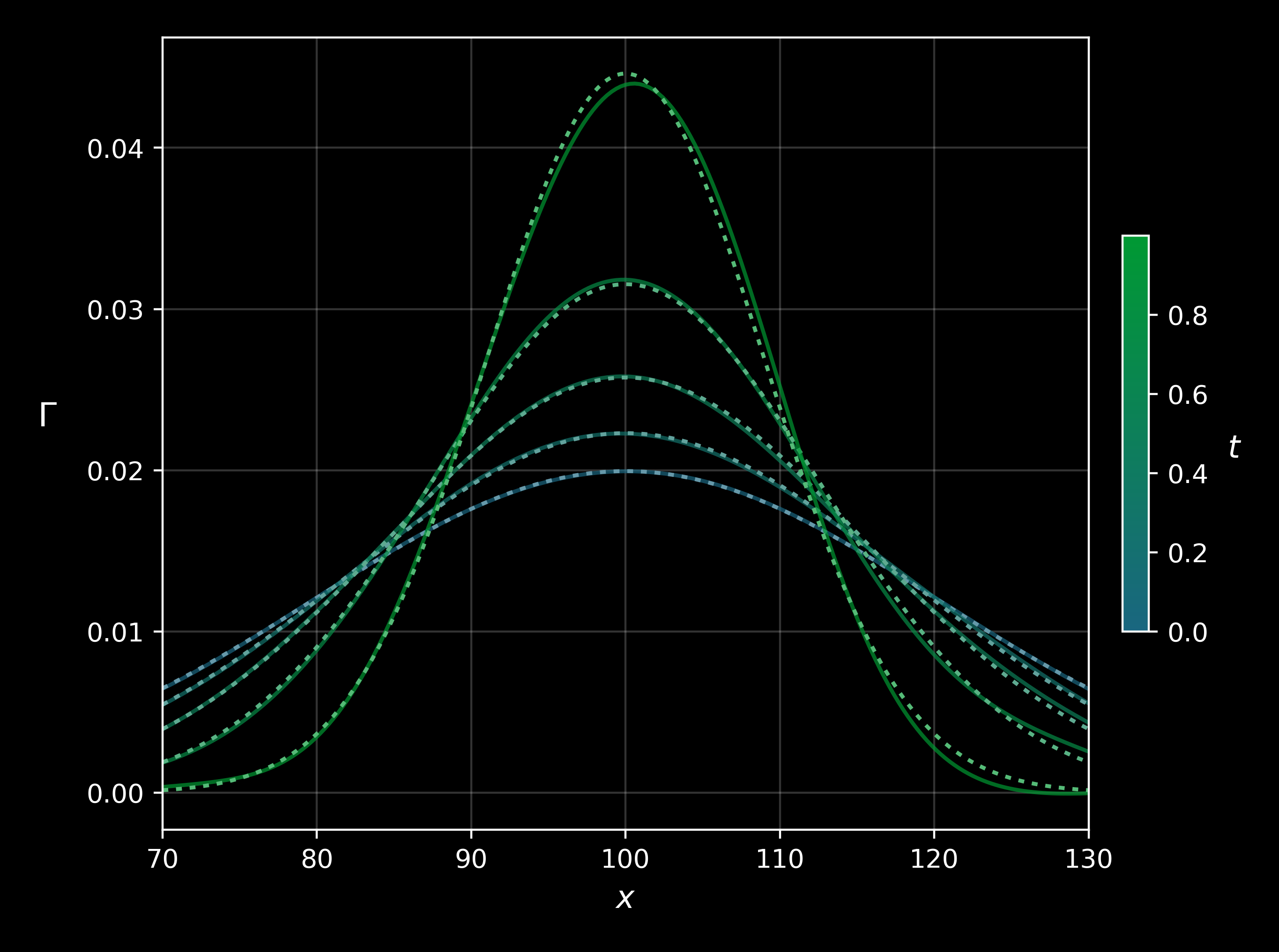}
        \label{fig:BachelierGamma}
    \end{subfigure}
    \label{fig:Bacheliersol}
\end{figure}

 \paragraph{Sensitivities.} By Section~\ref{sec:Greeks}, our framework yields the sensitivities directly, provided the kernel is sufficiently smooth. The estimates of the Delta $\Delta=\partial_x u$ and of the Gamma $\Gamma\coloneqq\partial_{xx}u$ are given by \eqref{eq:derivatives}. Figure~\ref{fig:Bacheliersol} shows the results for $n_1=n_2=2^{10}$ collocation points. The Delta is accurate (Figure~\ref{fig:BachelierDelta}). The Gamma becomes harder to approximate near maturity (Figure~\ref{fig:BachelierGamma}), as the payoff kink at \(x_0\) leads the  Gamma to concentrate into a Dirac mass as \(t\to T\).

\subsection{Multi-Asset Black--Scholes PDE}\label{sec:maxcall}

Consider a call option on the maximum of $d\ge2$ underlying assets (max-call option), with payoff
\[
g(x)=\Big(\max_{i=1,\ldots,d}x_i-K\Big)^+,
\]
and strike price $K>0$. We use the Black--Scholes model with symmetric, positively correlated assets with dynamics
\(
dS_t^i=\sigma S_t^i\,dW_t^i,
\)
where $W=(W^1,\ldots,W^d)$ is a Brownian motion in $\R^d$ with pairwise instantaneous correlation $\rho>0$. Then the price $u(t,x)=\E^\Q[g(S_T)\mid S_t=x]$ solves
\[
\begin{aligned}
\partial_tu+\frac{\sigma^2}{2}
\sum_{i,j=1}^d\rho_{ij}x_ix_j\,\partial_{ij}u
&=0,
&&\textnormal{on }\calX_1=(0,T)\times(0,\infty)^d,\\
u(T,\cdot)&=g,
&&\textnormal{on }\calX_2=\{T\}\times(0,\infty)^d,
\end{aligned}
\]
where $\rho_{ii}=1$ and $\rho_{ij}=\rho$ for $i\neq j$.

We use the parabolic kernel \eqref{eq:parabolicKernel} with the backward transformation $\psi(t,x)=\frac{x-S_0\mathbf 1}{S_0\sigma\sqrt{T-t+\varepsilon}}$, where $\mathbf 1=(1,\ldots,1)^\top$. Compared with \eqref{eq:transformBwd}, the price coordinates are centered at $S_0$ and standardized by the diffusion scale $S_0\sigma$. We smooth the terminal payoff with the log-sum-exp function:
\[
\operatorname{log\text{-}sum\text{-}exp}^{\varepsilon}(y_1,\ldots,y_n)
\coloneqq
\varepsilon\ln\Big(\sum_{i=1}^n e^{y_i/\varepsilon}\Big)
\approx \max(y_1,\ldots,y_n).
\]
Since $g(x) = \max \big(x_1 - K, \ldots ,x_d - K, 0\big)$, we apply the log-sum-exp function with $n=d+1$, $y_i = x_i - K$, and $y_{d+1} = 0$. This gives
\[
g(x)\approx g_\varepsilon(x)
\coloneqq
\varepsilon\ln\Big(1+\sum_{i=1}^d e^{(x_i-K)/\varepsilon}\Big).
\]
For the collocation points, we use Sobol sequences, see, e.g., \cite[Chapter~5]{Glasserman2004Monte}. Sobol points cover the domain and its boundary more evenly than random points. This typically improves the conditioning of the matrix $\bm{A}\bm{A}^*$. Regularization is still needed for a stable numerical solution. We map scrambled Sobol points to the uniform distribution on $(0,T)$ in time, and to the correlated lognormal distribution of $S_0e^{\sigma\sqrt T Z}$, with $Z\sim\calN(0,(\rho_{ij}))$, in the spatial variables. The terminal collocation points follow the same spatial distribution. The Sobol coordinates are clipped to $[10^{-6},1-10^{-6}]$ before applying the inverse normal distribution function, so that the spatial collocation points lie in a bounded set.

\begin{figure}
    \centering
    \caption{Max-call option on $d = 2$ assets at time $t=0$. Left: solution as a function of $(x_1,x_2)$. Right: relative error with respect to the Monte Carlo benchmark.}
    \includegraphics[width=0.7\linewidth]{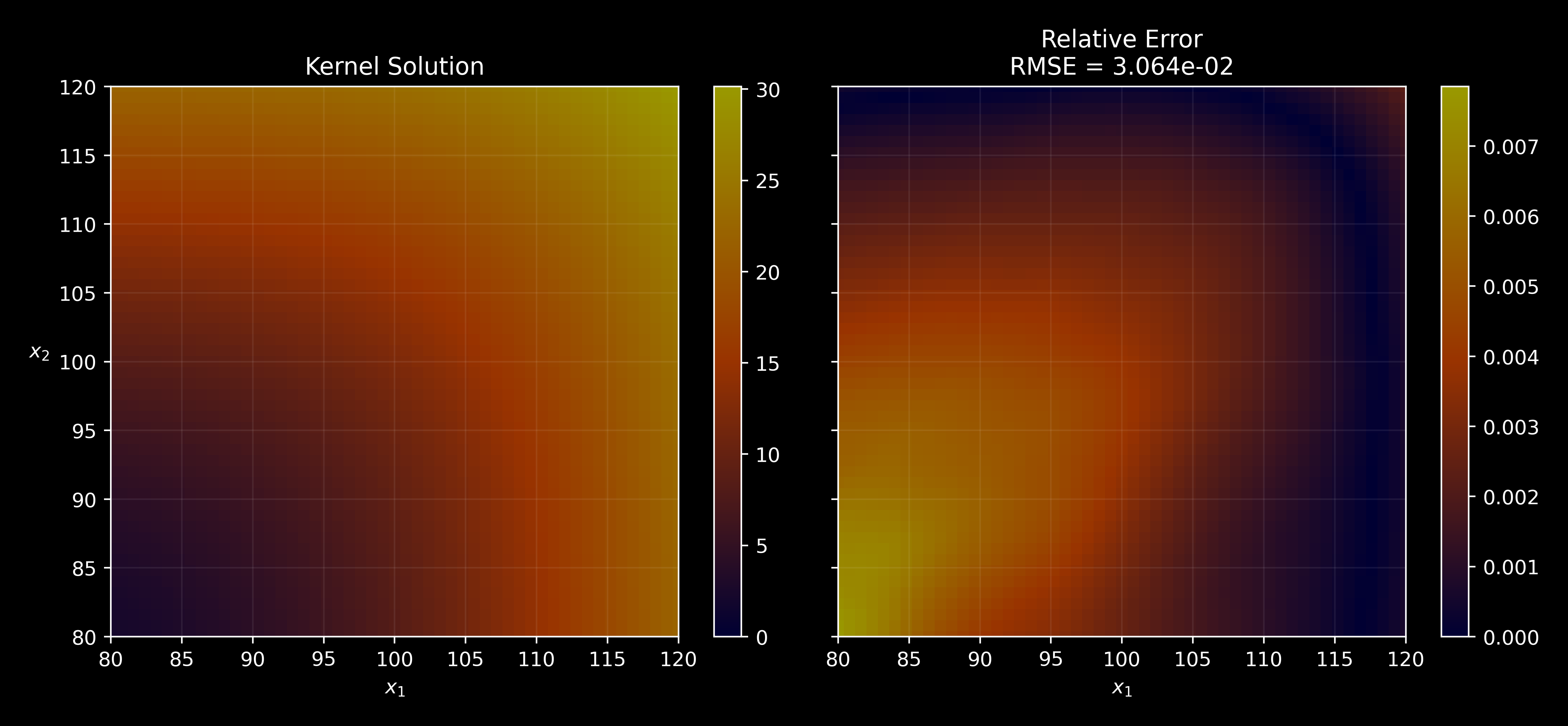}

    \label{fig:maxcallheatmap}
\end{figure}

\begin{figure}[H]
\caption{Max-call option on $d \in \{2,3,5\}$ assets. Top row: solution as a function of time $t$ at $x_1=\cdots=x_d=100$, compared with  Monte Carlo  (error bars). Bottom row: objective $\widehat\Phi_\lambda(\widehat h_{\lambda,r})$ of the reduced solution as a function of the number $r$ of selected representers.}
    \centering
      \begin{subfigure}[t]{0.325\linewidth}
        \centering
        \caption{$d = 2$}
        \includegraphics[width=0.99\textwidth]{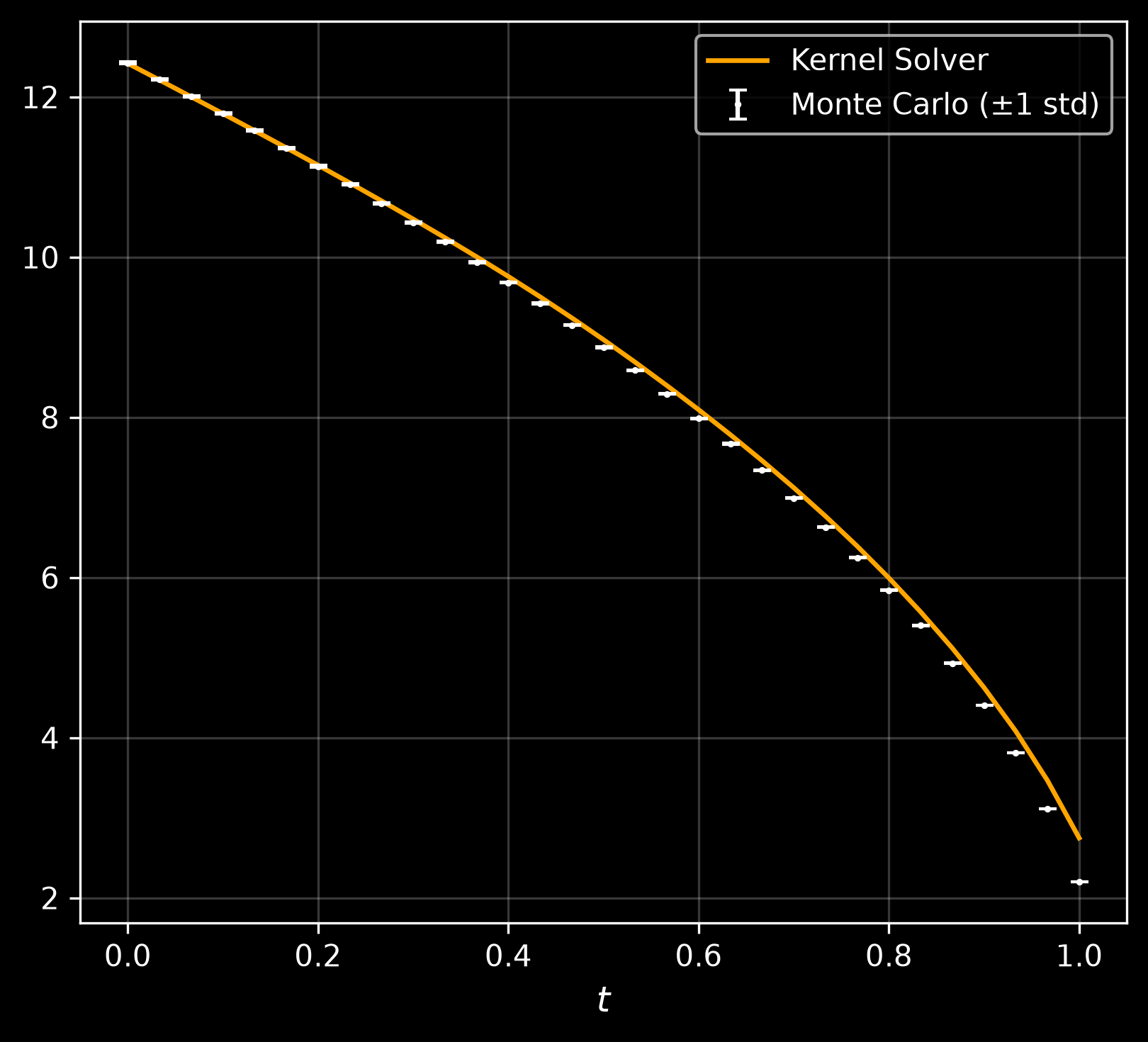}
        \label{fig:maxcall2}
    \end{subfigure}
    \begin{subfigure}[t]{0.325\linewidth}
        \centering
       \caption{$d = 3$}
        \includegraphics[width=0.99\textwidth]{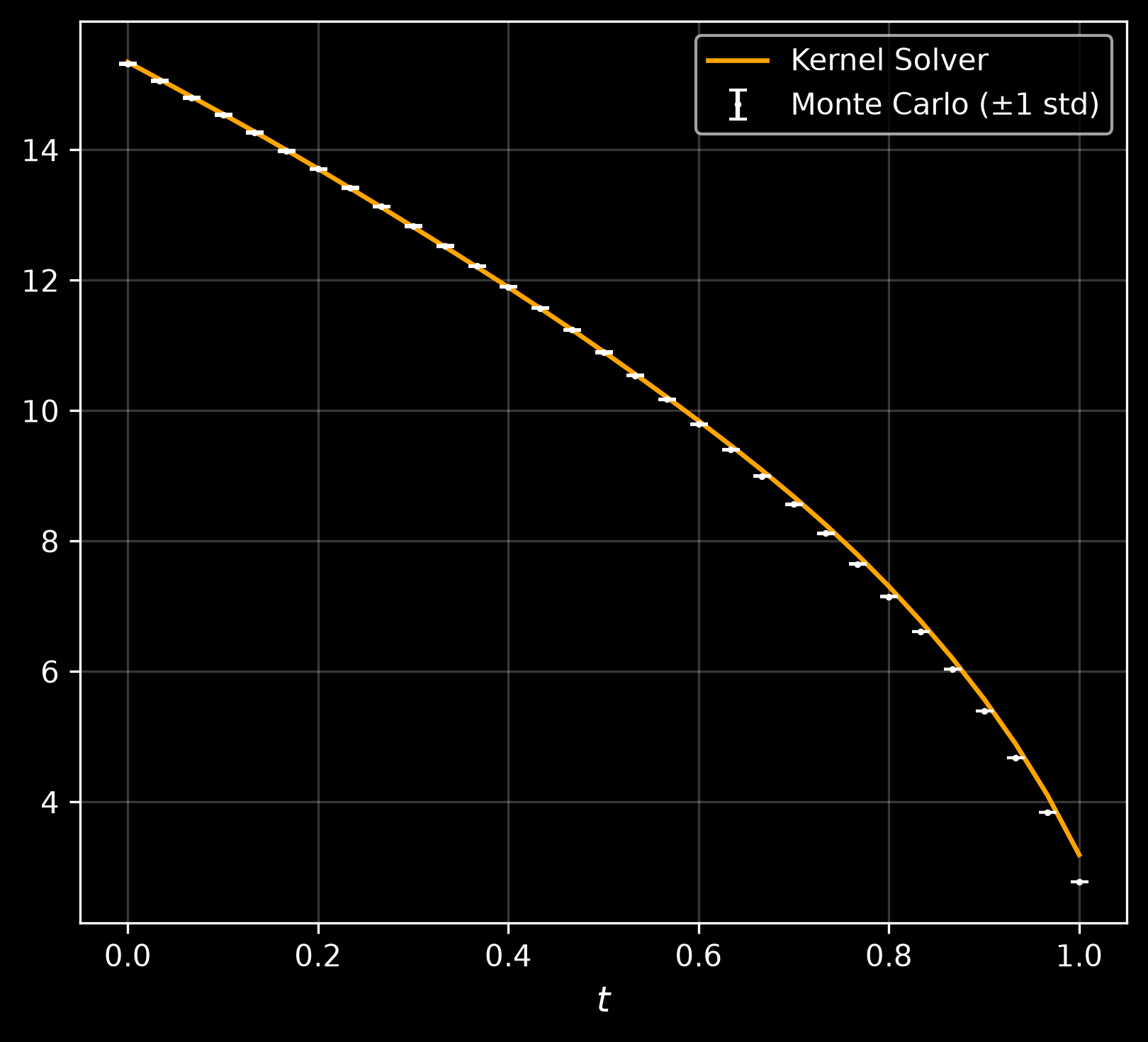}
        \label{fig:maxcall3}
    \end{subfigure}
        \begin{subfigure}[t]{0.325\linewidth}
        \centering
        \caption{$d = 5$}
        \includegraphics[width=0.99\textwidth]{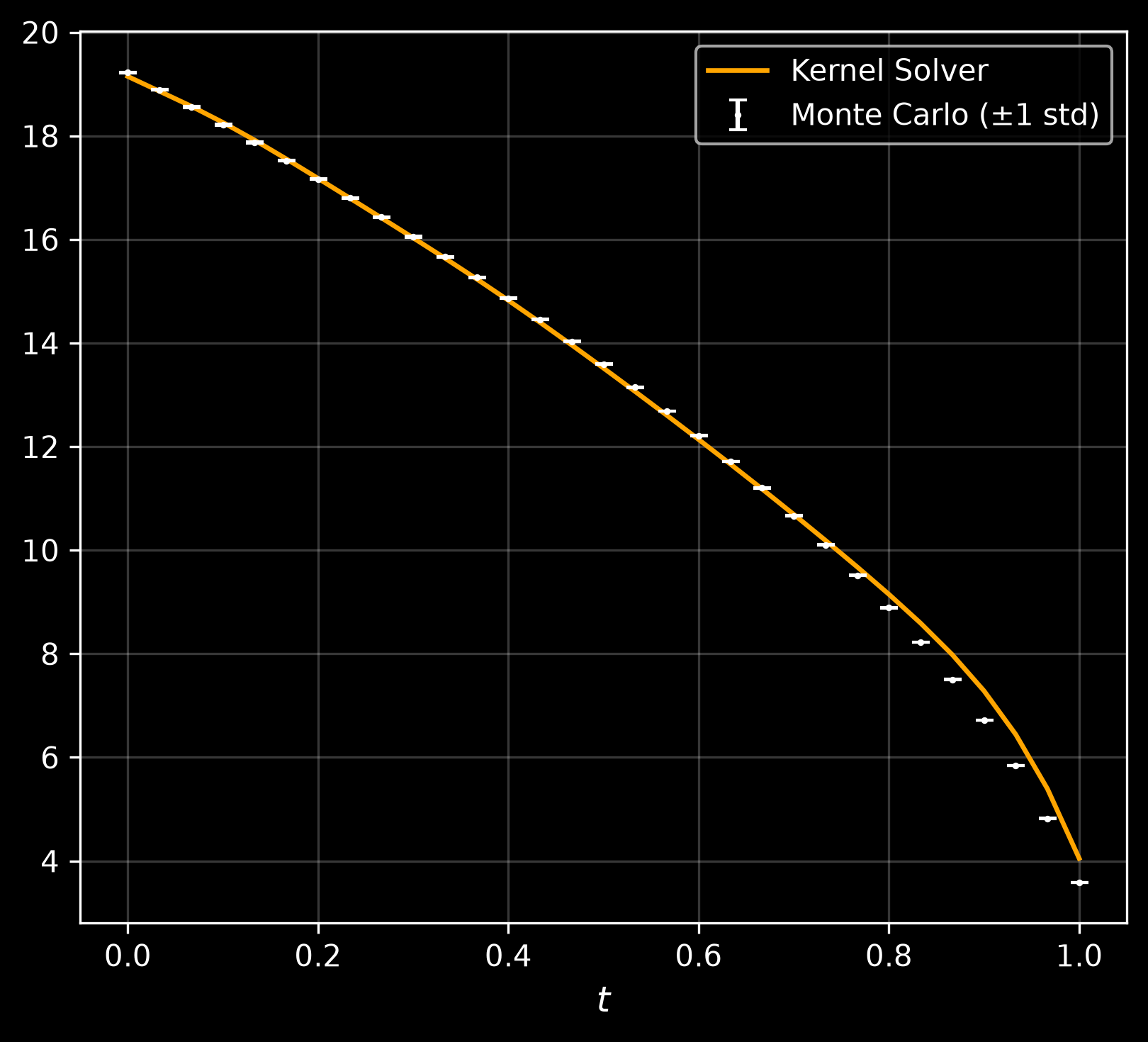}
        \label{fig:maxcall5}
    \end{subfigure}
    \label{fig:maxcallTimevalue}

    \vspace{-0.5em}

          \begin{subfigure}[t]{0.325\linewidth}
        \centering
        \caption{$d = 2$}
        \includegraphics[width=0.99\textwidth]{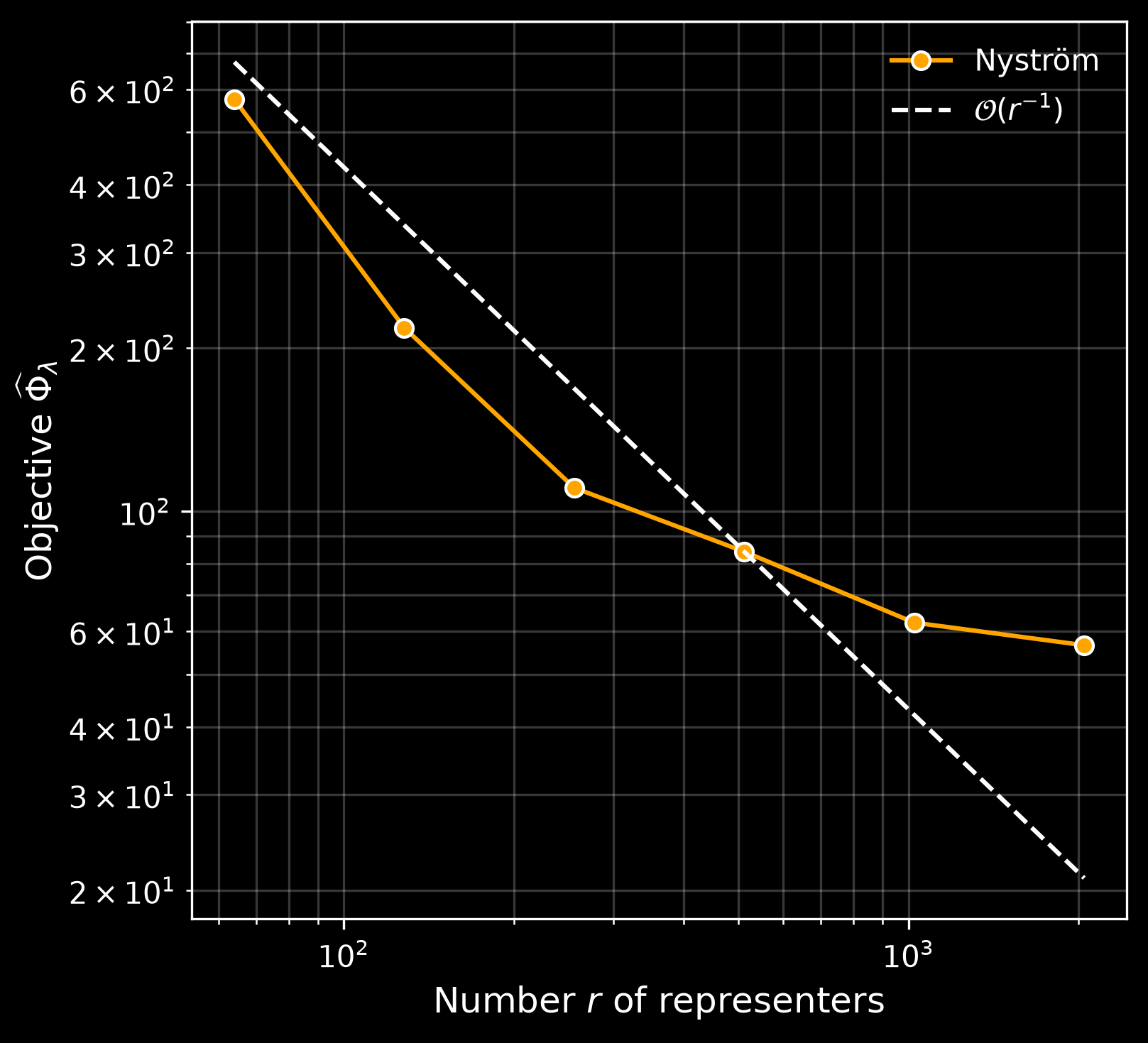}
        \label{fig:maxcallObj2}
    \end{subfigure}
    \begin{subfigure}[t]{0.325\linewidth}
        \centering
       \caption{$d = 3$}
        \includegraphics[width=0.99\textwidth]{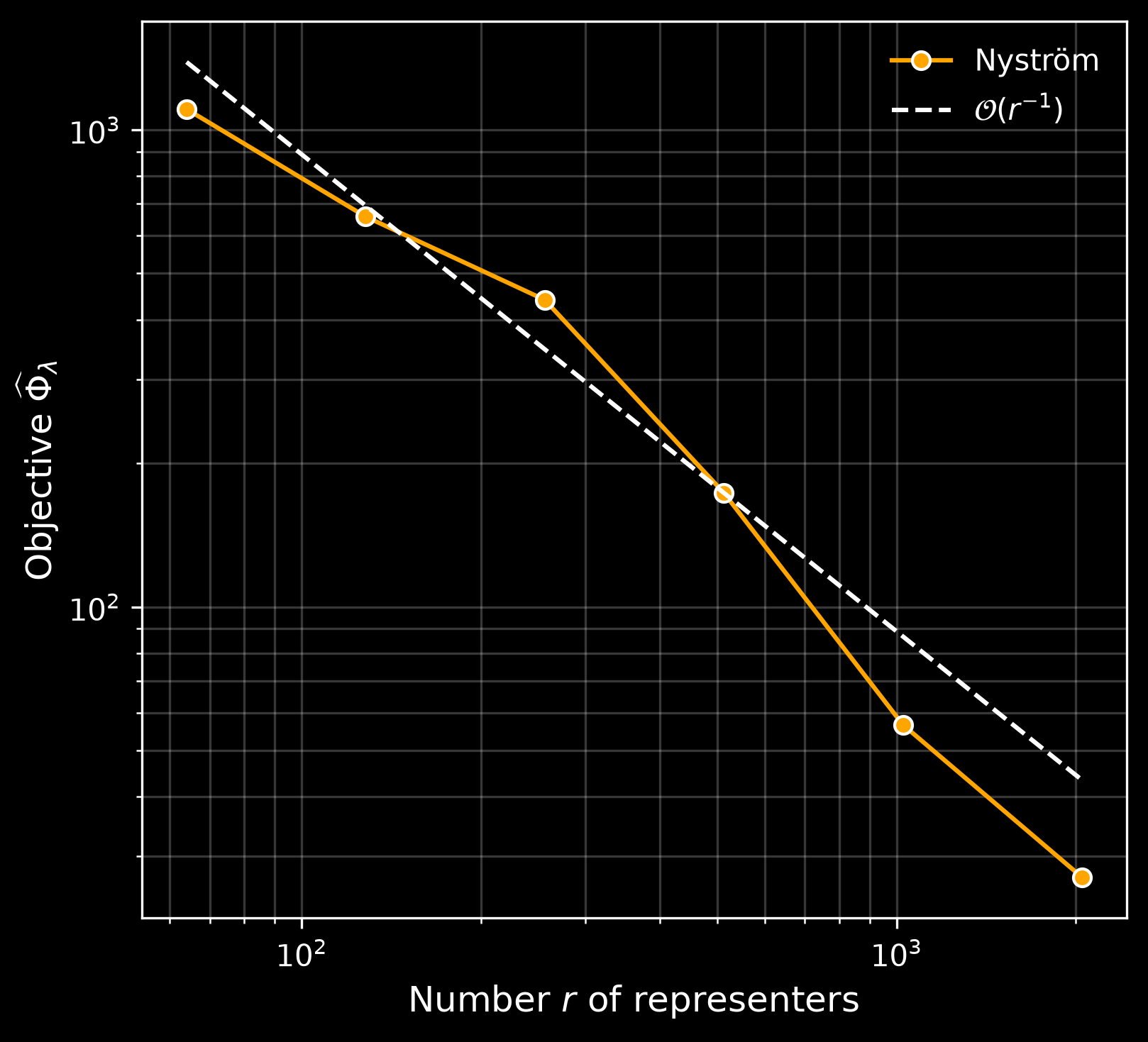}
        \label{fig:maxcallObj3}
    \end{subfigure}
        \begin{subfigure}[t]{0.325\linewidth}
        \centering
        \caption{$d = 5$}
        \includegraphics[width=0.99\textwidth]{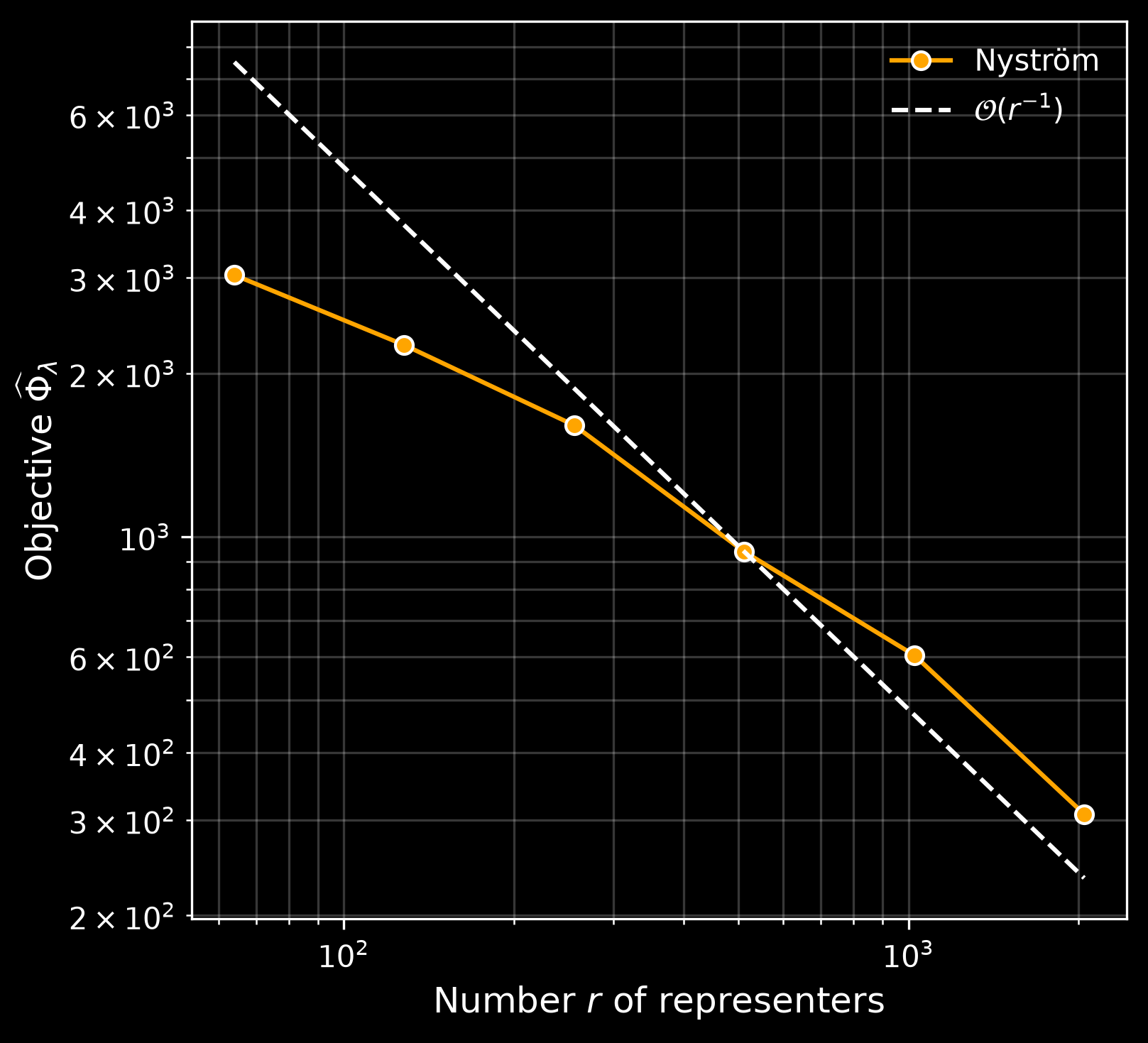}
        \label{fig:maxcallObj5}
    \end{subfigure}
\end{figure}

For all $d\in\{2,3,5\}$, we use the payoff and model parameters
\[
S_0 = K = 100, \qquad T=1, \qquad \sigma=0.20, \qquad \rho=0.50,
\]
\[
\varepsilon_{\mathrm{payoff}}=2, \quad
\varepsilon_{\mathrm{warp}}=0.10, \quad
\gamma_0=4, \quad
\gamma_1=\frac{1}{16}, \quad
w_1=1, \quad
w_2=10, \quad
\lambda=10^{-7}.
\]
Here $\varepsilon_{\mathrm{payoff}}$ denotes the smoothing parameter $\varepsilon$ in $g_\varepsilon$, and $\varepsilon_{\mathrm{warp}}$ the parameter $\varepsilon$ in the transformation $\psi$. We use $n_1=n_2=2^{13}$ collocation points and Nystr\"om approximations with
$r\in\{2^6,\ldots,2^{11}\}$ selected representers. Of these, $r_2=\lfloor 0.4\,r\rfloor$ are selected from the terminal condition and $r_1=r-r_2$ from the interior equation. The index sets are drawn uniformly at random, once. They are nested in $r$, being prefixes of one random permutation of each set of collocation points. All equations are formulated in price coordinates $x\in(0,\infty)^d$. The top row of Figure~\ref{fig:maxcallTimevalue} shows the reduced solution $\widehat h_{\lambda,r}$ with $r=2^{11}$. The Monte Carlo benchmark uses $2^{20}$ paths with antithetic variates. It is computed for the smoothed payoff $g_\varepsilon$, which is the terminal condition imposed on the kernel solution. At maturity it equals $g_\varepsilon(K,\ldots,K)=\varepsilon\ln(1+d)$, that is, $2.20$ for $d=2$ and $3.58$ for $d=5$. The deviation of the reduced solution from the benchmark is largest close to maturity. There, the spatial length scale of the kernel in price units, $\gamma_1^{-1/2}S_0\sigma\sqrt{\varepsilon_{\mathrm{warp}}}\approx 25$, remains large compared with the smoothing scale $\varepsilon_{\mathrm{payoff}}=2$ of the terminal payoff.

Figure~\ref{fig:maxcallheatmap} shows the kernel solution at time zero as a function of $(x_1,x_2)$, and its relative error with respect to a Monte Carlo benchmark. The top row of Figure~\ref{fig:maxcallTimevalue} compares the solution at the money, as a function of time $t$, with the Monte Carlo benchmark. The bottom row shows the objective $\widehat\Phi_{\lambda}(\widehat h_{\lambda,r})$ of the reduced solution as a function of the number $r$ of selected representers, together with an $r^{-1}$ reference line. Note that \cref{prop:NystromRate} bounds the expected objective gap $\widehat\Phi_\lambda(\widehat h_{\lambda,r})-\widehat\Phi_\lambda(\widehat h_\lambda)$. The panels show the objective itself, for the single draw of nested index sets described above. In Figure~\ref{fig:maxcallheatmap}, the benchmark is a Monte Carlo estimate with $2^{18}$ common paths for the original payoff $g$. The right panel shows the relative error $|\widehat u-u^{\mathrm{MC}}|/u^{\mathrm{MC}}$ on $[80,120]^2$, and the RMSE in its title is that of the absolute price error $\widehat u-u^{\mathrm{MC}}$ over the same grid. 

\subsection{Pricing PIDE in the Bates Model}\label{sec:jumpdiffusion}

This section shows how the kernel solver handles the integral operators arising in derivative-pricing equations under jump-diffusion dynamics. We consider the Bates model \cite{bates1996jumps}. 
Let $X_t = \ln S_t$ be the log-price and $V_t$ be the associated instantaneous variance. Assuming zero interest rates and dividends, the dynamics of $(X,V)$ in the Bates model under a risk-neutral measure $\Q$ are given by
\begin{align}
    dX_t &= -\big(V_t/2  + \zeta \gamma^{\Q} \big) dt + \sqrt{V_t} dW_t +  d \sum_{i=1}^{N_t} J_i, \\
    dV_t &= \kappa (\bar{V} - V_t) dt + \eta \sqrt{V_t} dZ_t,
\end{align}
where $W$ and $Z$ are Brownian motions with instantaneous correlation $\rho \in [-1,1]$, $N$ is an independent Poisson process with intensity $\zeta >0$, and $J_1,J_2,\ldots \overset{i.i.d.}{\sim} \mathcal{N}(\alpha, \beta)$ are the jump sizes of the log-price. The drift correction $\zeta\gamma^{\Q}$, with $\gamma^{\Q} = e^{\alpha + \frac{1}{2}\beta} - 1$, ensures that the asset price is a martingale under $\Q$.

The price of a European-style derivative with payoff $g: \R \to \R$ in this model satisfies the~PIDE
\begin{equation}
\partial_t u +  \calG u + \calI u - \zeta u = 0, \qquad u(T,x,v) = g(x),
\end{equation}
with the differential operator
\begin{equation}
\calG u={} -\big(v/2+\zeta\gamma^\Q\big)\partial_xu
+\kappa(\bar V-v)\partial_vu+\frac12v\,\partial_{xx}u
+\rho\eta v\,\partial_{xv}u
+\frac12\eta^2v\,\partial_{vv}u,
\end{equation}
and the integral operator
\begin{align}
\calI u(t,x,v)
&=\zeta\int_\R u(t,x+z,v)\,\phi_{\alpha,\beta}(z)\,dz,
\label{eq:jump_integral}\\[0.5em]
\phi_{\alpha,\beta}(z)
&=\frac{1}{\sqrt{2\pi\beta}}
\exp\!\left(-\frac{(z-\alpha)^2}{2\beta}\right).
\label{eq:GaussianDensity}
\end{align}

Since $\calI u(t,x,v)-\zeta u(t,x,v)=\zeta \int_{\R}\big(u(t,x+z,v)-u(t,x,v)\big) \phi_{\alpha,\beta}(z)\,dz$ is of the form \eqref{eq:jumpOperator} with a finite measure that acts on the coordinate $x$ only, \cref{prop:jump} applies to the Gaussian kernels used below.

For classical methods such as finite differences, the nonlocal term in \eqref{eq:jump_integral} adds computational cost, since it typically leads to dense matrices. In our framework, the action of the integral operator $\calI$ on Gaussian kernels is available in closed form. No separate discretization of the integral operator is needed. To see this, consider the spatial Gaussian kernel $k(x,y)=e^{-\frac{\gamma}{2}|x-y|^2}$, for fixed time $t$ and variance $v$. By the convention \eqref{kcalAdef}, the operator $\calI$ acts on the first argument of the kernel. Hence, with $c=\sqrt{\frac{2\pi}{\gamma}}$ and the Gaussian density $\phi$ in \eqref{eq:GaussianDensity},
\[
k_x^{\calI}(y)
=\calI k_y(x)
=\zeta\int_\R k(x+z,y)\phi_{\alpha,\beta}(z)\,dz
=\zeta c\,(\phi_{x,1/\gamma}*\phi_{\alpha,\beta})(y)
=\zeta c\,\phi_{x+\alpha,\beta+1/\gamma}(y).
\]
The cross-kernels in the blocks \eqref{eq:blocks}, such as $(k^{\calI}_x,k^{\calI}_y)_\calH$ and $(k^{\calI}_x,k^{\calG}_y)_\calH$, can likewise be computed explicitly. The same holds for the anisotropic Gaussian kernel
\[
k\big((t,x,v),(s,y,v')\big)
=
\exp\left\{
-\frac{\gamma_0}{2}|t-s|^2
-\frac{\gamma_1}{2}|x-y|^2
-\frac{\gamma_2}{2}|v-v'|^2
\right\}.
\]
We use this Gaussian kernel without the time-dependent transformation \eqref{eq:transformBwd}. Since the variance process is mean-reverting, no rescaling is needed in that direction. A time-dependent rescaling of the log-price coordinate did not improve the accuracy in our experiments.

\begin{figure}
    \centering
        \caption{Kernel price of a call option in the Bates model as a function of the underlying price $S=e^x$ and the variance $v$, at $t\in\{0,0.5,1\}$, using $n_1=n_2=2^{11}$ collocation points.}
    \includegraphics[width=0.95\linewidth]{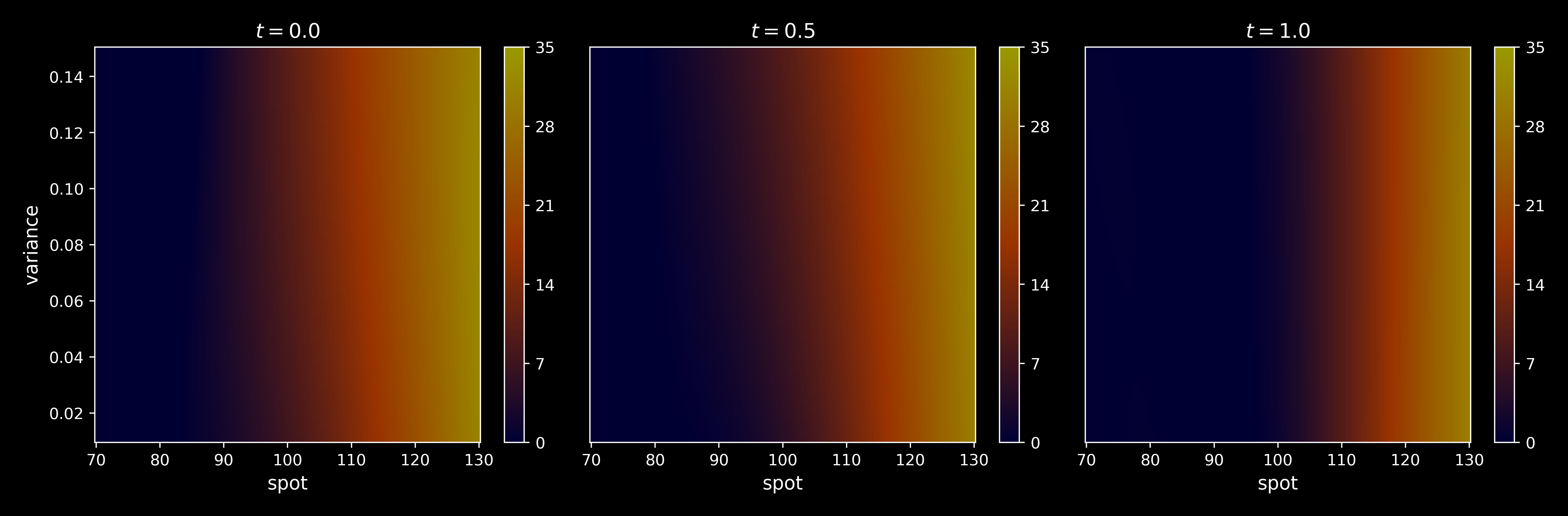}
    \label{fig:batesHeatmap}
\end{figure}

\begin{figure}
    \centering
     \caption{Kernel price of a one-year call option in the Bates model as a function of the underlying price $S = e^x$, at $t=0$ and $v=0.04$, compared with the FFT benchmark.}
    \includegraphics[width=0.5\linewidth]{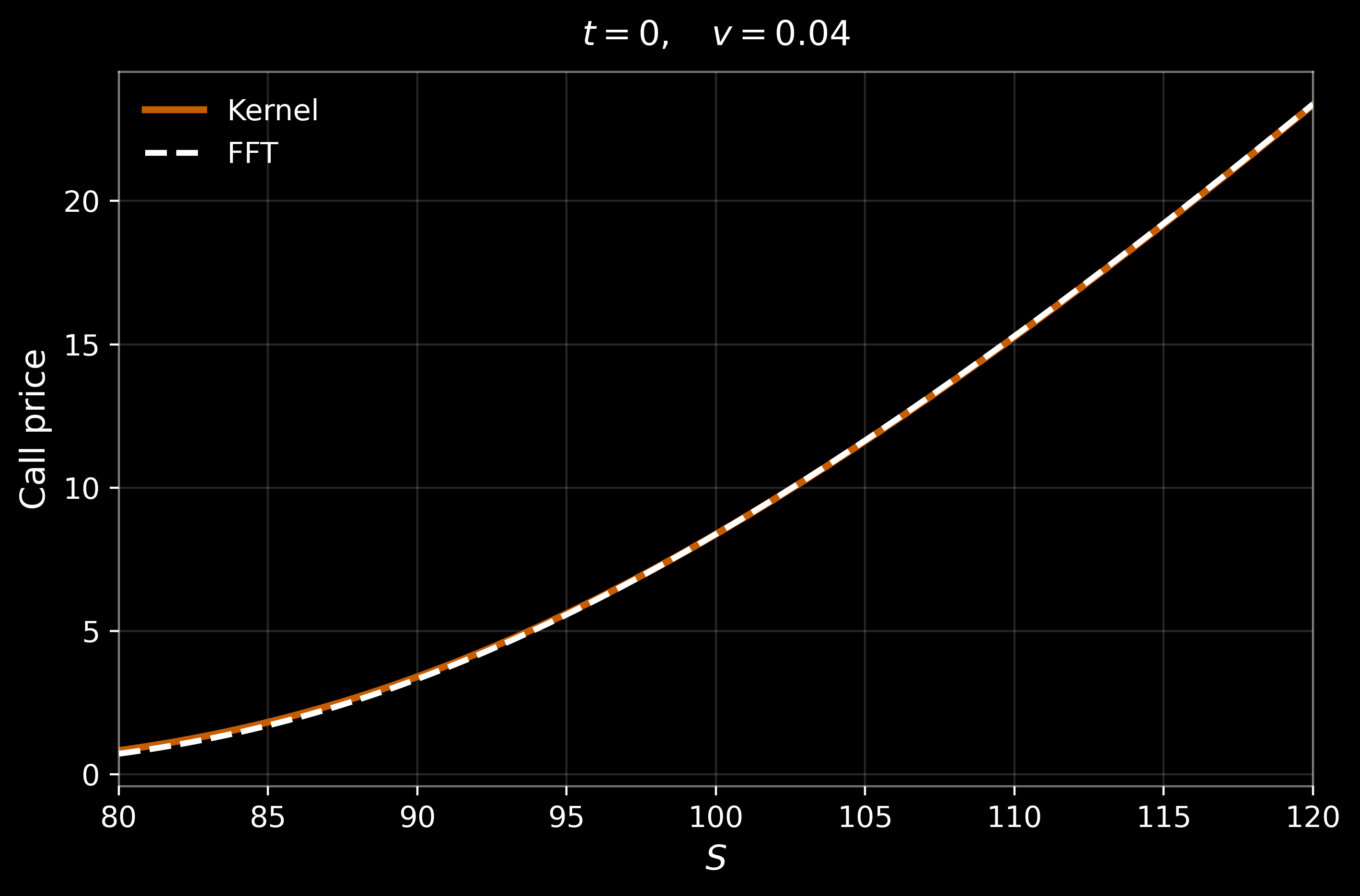}
    \label{fig:BatesSample}
\end{figure}

Figure~\ref{fig:batesHeatmap} shows the kernel solution, which we call the \emph{kernel price}, of a one-year call option with payoff $g(x)=(e^x-100)^+$ at three points in time. We use $n_1 = n_2 = 2^{11}$ scrambled Sobol collocation points each for the interior equation and the terminal condition, in the box $\ln S\in[\ln 100-0.6,\,\ln 100+0.6]$, $v\in[0.01,0.16]$. We use $\lambda=10^{-6}$ and $w_i =n_i$, so that $\lambda\bm{G}=\lambda I_N$. The model parameters are
\[
\zeta = 0.3, \quad
\alpha = -0.1, \quad
\beta = (0.15)^2, \quad
\kappa = 1.5, \quad
\bar V = 0.04, \quad
\eta = 0.3, \quad
\rho = -0.75,
\]
while the kernel parameters, selected by cross-validation over a grid of $(\gamma_0,\gamma_1,\gamma_2)$, are
$\gamma_0 = 1$, $\gamma_1 = 25$, and $\gamma_2 \approx 14$. 
Figure~\ref{fig:BatesSample} shows the kernel price of the one-year call option at $t=0$ as a function of the underlying price $S = e^x$, compared with the fast Fourier transform (FFT) benchmark of Carr and Madan~\cite{carrmadan}.

\subsection{Pricing PIDE in the Variance Gamma Model}\label{sec:VG}

Finally, we consider a pricing PIDE with a non-Gaussian L\'evy measure. The log-price follows, up to a deterministic drift, an asymmetric Variance Gamma (VG) process \cite{varianceGamma} with parameters $(\theta,\sigma, \nu)$. This is a pure-jump process of infinite activity with L\'evy measure
\begin{equation}\label{eq:levyVG}
\mu(dz)=\frac{e^{-\zeta(z)|z|}}{\nu|z|}\,dz,\qquad z\ne0, \quad 
\zeta(z)=\left(
\sqrt{\frac{\theta^2\nu^2}{4}+\frac{\sigma^2\nu}{2}}
+\frac{\theta\nu}{2}\operatorname{sign}(z)
\right)^{-1}.
\end{equation}
The price of a European-style derivative with payoff $g:\R\to \R$ solves the backward PIDE
\begin{equation}\label{eq:PIDEVG}
\begin{aligned}
\partial_tu+(r+\omega)\partial_xu
+\int_\R\big(u(t,x+z)-u(t,x)\big)\,\mu(dz)-ru&=0,\\
u(T,x)&=g(x).
\end{aligned}
\end{equation}
with the interest rate $r\ge0$ and the convexity correction $\omega = \frac{1}{\nu} \ln \big( 1 - \theta \nu - \frac{\sigma^2 \nu}{2} \big)$, where we assume that $\theta \nu + \frac{\sigma^2 \nu}{2}<1$.
As in \cref{sec:jumpdiffusion}, the variable $x$ represents the log-price of the underlying asset. The integral operator in \eqref{eq:PIDEVG} is of the form \eqref{eq:jumpOperator}, with jumps in the coordinate $x$ only, and the L\'evy measure \eqref{eq:levyVG} satisfies \eqref{eq:finiteVariation}. We use the Gaussian kernel
\[
k\big((t,x),(s,y)\big)
=\exp\left\{
-\frac{\gamma_0}{2}|t-s|^2
-\frac{\gamma_1}{2}|x-y|^2
\right\}
\]
on $\calS=[0,T]\times\R$, without the transformation $\psi$ in \eqref{eq:transformBwd}. Indeed, recall that $\psi$ is designed to match the scaling of diffusions, which is not natural for the pure-jump dynamics considered here. For this choice of kernel, \cref{prop:jump} applies, and \eqref{eq:jumpAssumption} holds with $c_0=1$ and $c_1=\sqrt{\gamma_1}$. In contrast to Section~\ref{sec:jumpdiffusion}, the integral operator applied to Gaussian kernels has no closed-form expression. The L\'evy measure has infinite activity, but it decays exponentially. We therefore approximate the integral in \eqref{eq:PIDEVG} by Gauss--Laguerre quadrature, applied separately on the positive and negative half-lines. No truncation of small jumps is needed.

\begin{figure}
    \centering
        \caption{Kernel price of a one-year, at-the-money call option in the VG model as a function of the common number $n=n_1=n_2$ of collocation points. Shaded area: one standard error across $16$ runs. Dotted line: FFT benchmark.}
    \includegraphics[width=0.48\linewidth]{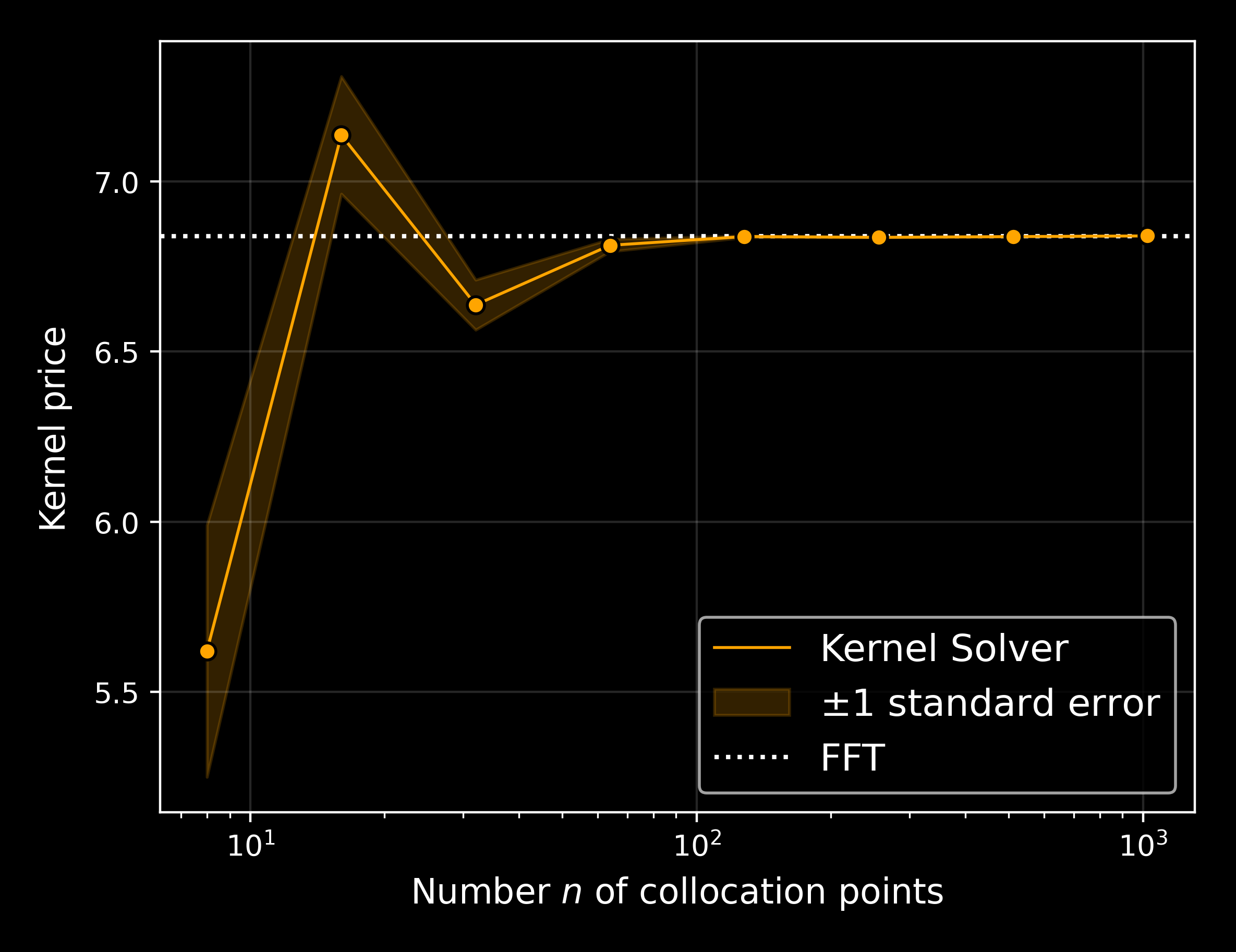}
    \label{fig:VGSample}
\end{figure}

Figure~\ref{fig:VGSample} shows the kernel price of a one-year, at-the-money call option in the VG model as a function of the common number $n=n_1=n_2$ of interior and terminal collocation points. The shaded area shows the mean plus and minus one standard error across $16$ runs of the kernel solver with scrambled Sobol sequences. For large $n$, the kernel price stabilizes and approaches the FFT benchmark of \cite{carrmadan}. It also becomes less sensitive to the choice of the scrambled Sobol points. We use the VG parameters $\theta=-0.14$, $\sigma=0.12$, $\nu=0.2$, the interest rate $r=0.03$,  $S_0=K=100$ (at-the-money), and $T=1$. The collocation points lie in the interval $[\ln S_0-5\sigma\sqrt T,\,\ln S_0+5\sigma\sqrt T]$. The kernel parameters are $\gamma_0=4$ and $\gamma_1=(2\sigma\sqrt T)^{-2}\approx17.4$. We use $\lambda=10^{-7}$ and $w_1=w_2=n$, so that $\lambda\bm{G}=\lambda I_N$. The quadrature uses $12$ Gauss--Laguerre nodes on each half-line.

Table~\ref{tab:VG} reports a $99\%$ confidence interval (CI) for the mean kernel price across $16$ randomized runs, each with $n=2^{10}$ collocation points. The interval is narrow and contains the FFT benchmark.

\begin{table}
    \centering
        \caption{$99\%$ confidence interval (CI) for the mean kernel price of a one-year, at-the-money call option in the VG model across $16$ randomized runs with $n_1=n_2=2^{10}$ collocation points, compared with the FFT benchmark of \cite{carrmadan}.}
    \begin{tabular}{c|c}
       Kernel ($99\%$ CI) & FFT benchmark \\ \hline
       $[6.8395, 6.8403]$  & $6.8398$
    \end{tabular}
    \label{tab:VG}
\end{table}

\section{Conclusion}\label{sec:conclusion}

We developed a mesh-free kernel solver for systems of linear operator equations, treating differential, boundary, and integral constraints within a common regularized least-squares framework. Sensitivities follow by differentiation of the kernel representation. To reduce computational cost, we introduced a Nystr\"om approximation based on subsampling the operator-induced representers. The reduced solution satisfies a best-approximation property in the regularized empirical norm, with an expected objective gap of order $r^{-1}$ under balanced uniform subsampling. Numerical experiments on multidimensional PDEs and pricing PIDEs with jump dynamics support these findings and illustrate the accuracy and flexibility of the approach. 
Future work includes scaling the framework to higher-dimensional problems, where further randomized approximations and GPU implementations may be needed to reduce the cost of dense kernel computations.

\appendix

\section{Proofs}\label{app:proofs}

\subsection{Proofs for \cref{sec:prelim}}

\begin{proof}[Proof of \cref{lem:ARegular}]
Let $\calH$ be $\calA$-regular and fix $x\in\calX$. By the Riesz representation theorem, there exists $q_x\in\calH$ such that $\calA h(x)=(q_x,h)_\calH$ for all $h\in\calH$. Taking $h=k_y$ and using \eqref{eq:reproduce} gives $q_x(y)=(q_x,k_y)_\calH=\calA k_y(x)=k^\calA_x(y)$ for all $y\in\calS$. Hence $k^\calA_x=q_x\in\calH$ and \eqref{reprAh} holds. Conversely, \eqref{reprAh} and the Cauchy--Schwarz inequality give $|\calA h(x)|\le\|k^\calA_x\|_\calH\|h\|_\calH$, so that $h\mapsto\calA h(x)$ is continuous on $\calH$.
\end{proof}

\begin{proof}[Proof of \cref{propdiff}]
\ref{propdiff1}: Let $x \in \mathcal{X}$, write $v\in\R^d$ for the $i$-th standard basis vector, and choose $\varepsilon_0 >0$ such that $x + \varepsilon v \in \mathcal{X}$ for all $|\varepsilon|\le\varepsilon_0$. For $0<|\varepsilon|\le \varepsilon_0$, consider the difference quotients
    $q_\varepsilon \coloneqq \frac{1}{\varepsilon}(k_{x+\varepsilon v} - k_x) \in \calH.$
    By the reproducing property \eqref{eq:reproduce} and the fundamental theorem of calculus, applied first in the second and then in the first argument of $k$, we have for $0<|\varepsilon|,|\delta|\le\varepsilon_0$,
    \begin{align*}
(q_\varepsilon,q_\delta)_\calH
&=\frac{1}{\varepsilon\delta}\Big(
 k(x+\varepsilon v,x+\delta v)-k(x+\varepsilon v,x)\\
&\hspace{7em}-k(x,x+\delta v)+k(x,x)\Big)\\
&=\frac{1}{\varepsilon\delta}
\int_0^\varepsilon\int_0^\delta
\partial_{x_i}\partial_{y_i}k(x+sv,x+tv)\,dt\,ds.
    \end{align*}
    Since $\partial_{x_i}\partial_{y_i}k$ is continuous at $(x,x)$, the right-hand side converges to $\partial_{x_i}\partial_{y_i}k(x,x)$ as $\varepsilon,\delta\to0$. Hence
    \[
\|q_\varepsilon-q_\delta\|_\calH^2
=(q_\varepsilon,q_\varepsilon)_\calH
-2(q_\varepsilon,q_\delta)_\calH
+(q_\delta,q_\delta)_\calH
\longrightarrow0.
    \]
    Thus $q_\varepsilon$ converges in $\calH$, as $\varepsilon\to0$, to some $q\in\calH$ with $\|q\|_\calH^2=\partial_{x_i y_i}k(x,y)|_{y = x}$.

    For $h\in\calH$, the convergence in $\calH$ implies
    \begin{align*}
(q,h)_\calH
=\lim_{\varepsilon\to0}(q_\varepsilon,h)_\calH =\lim_{\varepsilon\to0}
\frac{h(x+\varepsilon v)-h(x)}{\varepsilon}.
    \end{align*}
    Hence the partial derivative $\partial_ih(x)$ exists and equals $(q,h)_\calH$. This shows that $\calH\subseteq\calD(\partial_i)$ and that $h\mapsto\partial_ih(x)$ is continuous on $\calH$, so that $\calH$ is $\partial_i$-regular. By \cref{lem:ARegular}, $q=k^{\partial_i}_x$, and $k^{\partial_i}_x(y)=\partial_ik_y(x)=\partial_{x_i}k(x,y)$ for all $y\in\calS$ by \eqref{kcalAdef}.

\ref{propdiff2}: By part \ref{propdiff1}, $\calH$ is $\partial_j$-regular. For $x,x'\in\calX$, the reproducing property \eqref{reprAh} of $k^{\partial_j}_{x'}$ and the symmetry of $k$ and of the inner product give
\[
(k^{\partial_j}_{x},k^{\partial_j}_{x'})_\calH
=\partial_jk^{\partial_j}_{x}(x')
=\widetilde k(x',x)
=\widetilde k(x,x').
\]
We can thus repeat the argument of part \ref{propdiff1} with the difference quotients $p_\varepsilon\coloneqq(k^{\partial_j}_{x+\varepsilon v}-k^{\partial_j}_x)/\varepsilon$ and with $\widetilde k$ in place of $k$. It follows that $p_\varepsilon$ converges in $\calH$ to some $p$ with $\|p\|_\calH^2=\partial_{x_iy_i}\widetilde k(x,y)|_{y = x}$, and that, for all $h\in\calH$,
\begin{align*}
(p,h)_\calH =\lim_{\varepsilon\to0}
\frac{\partial_jh(x+\varepsilon v)-\partial_jh(x)}{\varepsilon}
=\partial_{ij}h(x).
\end{align*}
Hence $\calH$ is $\partial_{ij}$-regular with $k^{\partial_{ij}}_x=p$, and $k^{\partial_{ij}}_x(y)=\partial_{ij}k_y(x)=\partial_{x_ix_j}k(x,y)$ for all $y\in\calS$.
\end{proof}

\begin{proof}[Proof of \cref{prop:jump}]
Fix $x\in\calX$. Since $k$ is continuous, the map $y\mapsto k_y:\calS\to\calH$ is continuous, because $\|k_y-k_{y'}\|_\calH^2=k(y,y)-2k(y,y')+k(y',y')$. Hence $z\mapsto k_{x+z}-k_x$ is continuous, and thus strongly measurable, on the Borel set $\{z:x+z\in\calS\}$. The complement of this set is a $\mu$-null set. By \eqref{eq:jumpAssumption},
$$\|k_{x+z}-k_x\|_\calH\le\min(2c_0,c_1\|z\|_2)$$
for $\mu$-almost all $z$, where we used that $\|k_y\|_\calH^2=k(y,y)\le c_0^2$. In view of \eqref{eq:finiteVariation}, the Bochner integral $q_x\coloneqq\int(k_{x+z}-k_x)\,\mu(dz)$ is a well-defined element of $\calH$, with norm bounded as in \eqref{eq:jumpRepresenter}. For $h\in\calH$, the reproducing property \eqref{eq:reproduce} gives 
\[|h(x+z)-h(x)|=|(k_{x+z}-k_x,h)_\calH|\le\min(2c_0,c_1\|z\|_2)\|h\|_\calH.\]
Hence $h\in\calD(\calA)$, and interchanging the inner product with the Bochner integral yields
$$\calA h(x)=\int(k_{x+z}-k_x,h)_\calH\,\mu(dz)=(q_x,h)_\calH.$$
Therefore $\calH$ is $\calA$-regular, and $q_x=k^\calA_x$ by \cref{lem:ARegular}.
\end{proof}

\begin{proof}[Proof of \cref{lem:L2}]
    \textit{1.} For $h$ in the linear span of $\{k_y:y\in\calS\}$, the function $x\mapsto(k^\calA_x,h)_\calH$ is a finite linear combination of the functions $x\mapsto k^\calA(x,y)$ and thus measurable. Since this span is dense in $\calH$, the function $x\mapsto(k^\calA_x,h)_\calH$ is measurable for every $h\in\calH$, as a pointwise limit of measurable functions. As $\calH$ is separable, Pettis' measurability theorem implies that $x\mapsto k^\calA_x$ is measurable. For $g\in L^2(\Pb)$, the Cauchy--Schwarz inequality and \eqref{assL2} give $\int_\calX\|k^\calA_x\|_\calH\,|g(x)|\,\Pb(dx)<\infty$, so that the Bochner integral is well defined.

    \textit{2.} For all $h\in \calH$, observe that $| \calA h(x)| = |( k^\calA_x,  h )_{\calH}|\le \| k^\calA_x \|_{\calH} \| h \|_{\calH} $, hence
    \begin{align*}
\|\calA h\|_{L^2(\Pb)}^2
=\int_\calX |(k^\calA_x,h)_\calH|^2\,\Pb(dx)\le\|h\|_\calH^2
\int_\calX\|k_x^\calA\|_\calH^2\,\Pb(dx)
<\infty.
    \end{align*}

    \textit{3.} Let $\Psi$ be a (countable) orthonormal basis of $\calH$. Then,
 the monotone convergence theorem and Parseval's identity imply  that
    \begin{align*}
\|\calA\|_{\textnormal{HS}}^2
=\sum_{\psi\in\Psi}\|\calA\psi\|_{L^2(\Pb)}^2=\int_\calX\sum_{\psi\in\Psi}
|(k^\calA_x,\psi)_\calH|^2\,\Pb(dx)=\int_\calX\|k^\calA_x\|_\calH^2\,\Pb(dx).
    \end{align*}

    \textit{4.} For all $h\in \calH$, $g\in L^2(\Pb)$, observe that
\begin{align*}
    ( \calA h, g)_{L^2(\Pb)} = \int_{\calX}   ( k^{\calA}_x, h)_{\calH} \ g(x) \ \Pb(dx)
    = \Big( h, \int_{\calX}    k^{\calA}_x\, g(x) \ \Pb(dx) \Big)_{\calH},
\end{align*}
which proves \eqref{eq:adjoint}.
\end{proof}

\subsection{Proofs for \cref{sec:linearOperatorEqs}}

The population and empirical problems are both instances of the following elementary result, which we state separately.

\begin{lemma}\label{lem:quadratic}
Let $E$ be a Hilbert space, $B:\calH\to E$ a bounded linear operator, $b\in E$, and $W$ a bounded, self-adjoint, positive semidefinite operator on $E$. For $\lambda>0$, consider the functional
\[
\Psi(h)\coloneqq
\big(W(Bh-b),Bh-b\big)_E
+\lambda\|h\|_\calH^2,
\qquad h\in\calH.
\]
Then $\Psi$ has a unique minimizer $h^\star$ over $\calH$, given by
\begin{equation}\label{eq:quadraticMin}
h^\star=(B^*WB+\lambda \bm I)^{-1}B^*Wb,
\end{equation}
and
\begin{equation}\label{eq:quadraticGap}
\Psi(h)-\Psi(h^\star)=\big((B^*WB+\lambda \bm I)(h-h^\star),h-h^\star\big)_\calH,\quad h\in\calH.
\end{equation}
If, in addition, $(We,e)_E\ge w_0\|e\|_E^2$ for all $e\in E$, for some $w_0>0$, then $W$ is invertible and
\begin{equation}\label{eq:quadraticMinDual}
h^\star=B^*(BB^*+\lambda W^{-1})^{-1}b.
\end{equation}
\end{lemma}

\begin{proof}
The operator $T\coloneqq B^*WB$ is bounded, self-adjoint, and positive semidefinite on $\calH$. Hence $T+\lambda\bm I$ is invertible, and we define $h^\star\coloneqq(T+\lambda \bm I)^{-1}B^*Wb$. Expanding the square and using $B^*Wb=(T+\lambda\bm I)h^\star$, we obtain
\begin{align*}
\Psi(h)
=\big((T+\lambda\bm I)h,h\big)_\calH
-2\big((T+\lambda\bm I)h^\star,h\big)_\calH+(Wb,b)_E,
\qquad h\in\calH,
\end{align*}
which implies \eqref{eq:quadraticGap}. Since the right-hand side of \eqref{eq:quadraticGap} is bounded from below by $\lambda\|h-h^\star\|_\calH^2$, $h^\star$ is the unique minimizer of $\Psi$. Under the additional assumption, $W$ has a bounded inverse $W^{-1}$, which is self-adjoint and bounded from below by $\|W\|^{-1}$. Hence $BB^*+\lambda W^{-1}$ is bounded, self-adjoint, and bounded from below by $\lambda\|W\|^{-1}$, thus invertible. Multiplying the identity
$(T+\lambda \bm I) B^* = B^*W (BB^* + \lambda W^{-1})$ from the left by $(T+\lambda\bm I)^{-1}$ and from the right by $(BB^* + \lambda W^{-1})^{-1}$ yields \eqref{eq:quadraticMinDual}.
\end{proof}

\begin{proof}[Proof of \cref{prop:population}]
By \cref{asm:L2} and \cref{lem:L2}, $\bm\calA$ is a bounded linear operator from $\calH$ to the Hilbert space $E=L^2(\Pb_1)\times\cdots\times L^2(\Pb_M)$, and $\Phi_\lambda$ is of the form considered in \cref{lem:quadratic} with $B=\bm\calA$, $b=\bm f$, and $W=\bm\calW$. Since $(\bm\calW\bm e,\bm e)_E=\sum_{m=1}^Mw_m\|e_m\|_{L^2(\Pb_m)}^2\ge\min_mw_m\,\|\bm e\|_E^2$ for $\bm e=(e_1,\ldots,e_M)\in E$, the additional assumption of the lemma holds with $w_0=\min_mw_m>0$, and \eqref{eq:hLambda}--\eqref{eq:hLambda2} follow from \eqref{eq:quadraticMin} and \eqref{eq:quadraticMinDual}.
\end{proof}

\begin{proof}[Proof of \cref{thm:representer}]
The sampling operator $\bm{A}:\calH\to\R^N$ is bounded, because $\calH$ is $\calA_m$-regular for all $m$. The functional $\widehat\Phi_\lambda$ is thus of the form considered in \cref{lem:quadratic} with $E=\R^N$, $B=\bm{A}$, $b=\bm y$, and the diagonal matrix $W=\bm{G}^{-1}$, whose diagonal entries $w_m/n_m$ are positive. Hence \eqref{eq:hLambda3} follows from \eqref{eq:quadraticMinDual}, and the representation \eqref{eq:representer} from \eqref{eq:hLambda3} and \eqref{eq:adjointEmp}.
\end{proof}

\begin{proof}[Proof of \cref{prop:Nystrom}]
As $\bm{A}^*\bm{G}^{-1}\bm{A}=\widehat T$, identity \eqref{eq:quadraticGap} in \cref{lem:quadratic} reads
\begin{equation}\label{eq:gapIdentity}
\widehat\Phi_\lambda(h)-\widehat\Phi_\lambda(\widehat h_\lambda)
=
\|h-\widehat h_\lambda\|_{\widehat T,\lambda}^2,
\qquad h\in\calH.
\end{equation}
Minimizing $\widehat\Phi_\lambda$ over the closed subspace $\calV$ is therefore equivalent to minimizing the distance to $\widehat h_\lambda$ in the norm $\|\cdot\|_{\widehat T,\lambda}$, which is equivalent to the norm of $\calH$ by \eqref{eq:normEquivalence}. The minimizer exists, is unique, and equals the orthogonal projection of $\widehat h_\lambda$ onto $\calV$ with respect to $(\cdot,\cdot)_{\widehat T,\lambda}$. This proves \eqref{eq:NystromBestApprox}, and \eqref{eq:NystromObjectiveGap} follows by taking $h=\widehat h_{\lambda,\calV}$ in \eqref{eq:gapIdentity}. Finally,
\eqref{eq:NystromHBound}--\eqref{eq:NystromResidualBound}
follow directly from \eqref{eq:empiricalNorm}.
\end{proof}

\begin{proof}[Proof of \cref{prop:NystromRate}]
By Proposition~\ref{prop:Nystrom},
$\widehat h_{\lambda,r}$ is the best approximation of
$\widehat h_\lambda$ in $\calH_r$ for the
$\|\cdot\|_{\widehat T,\lambda}$ norm. Hence, for every
random element $\widetilde h_r$ of $\calH_r$,
\begin{equation}\label{eq:NystromCandidate}
\widehat\Phi_\lambda(\widehat h_{\lambda,r})
-\widehat\Phi_\lambda(\widehat h_\lambda)
\le
\|\widetilde h_r-\widehat h_\lambda\|_{\widehat T,\lambda}^2.
\end{equation}
We choose the reweighted subsample
\(
\widetilde h_r
\coloneqq
\sum_{m=1}^M
\frac{n_m}{r_m}
\sum_{n\in\mathcal J_m}
\widehat\beta_{m,n}\varphi_{m,n}
\in\calH_r.
\)
For each $m$, set
\(
z_{m,n}
\coloneqq
n_m\widehat\beta_{m,n}\varphi_{m,n}\) and \(
\bar z_m
\coloneqq
\frac1{n_m}\sum_{n=1}^{n_m}z_{m,n}
=
\widehat h_{\lambda,m}.
\)
Then
\[
\widetilde h_r-\widehat h_\lambda
=
\sum_{m=1}^M
\left(
\frac1{r_m}\sum_{n\in\mathcal J_m}z_{m,n}
-\bar z_m
\right).
\]
Each term in the sum has zero $\E_{\mathcal J}$-expectation. Since the index sets $\mathcal J_m$ are independent across $m$, the cross terms vanish upon applying $\E_{\mathcal J}$. Therefore,
\[
\E_{\mathcal J}\!
\left[
\|\widetilde h_r-\widehat h_\lambda\|_{\widehat T,\lambda}^2
\right]
=
\sum_{m=1}^M
\E_{\mathcal J}\!
\left[
\left\|
\frac1{r_m}\sum_{n\in\mathcal J_m}z_{m,n}
-\bar z_m
\right\|_{\widehat T,\lambda}^2
\right].
\]

For a subset $J$ drawn uniformly at random among all subsets of $\{1,\ldots,p\}$ of size $q$, with $p\ge2$,
and arbitrary elements $z_1,\ldots,z_p$ of a Hilbert space with norm $\|\cdot\|$,
\begin{equation}\label{eq:finitePopulationVariance}
\E\left[
\left\|
\frac1q\sum_{i\in J}z_i-\bar z
\right\|^2
\right]
=
\frac{p-q}{q(p-1)}
\frac1p\sum_{i=1}^p\|z_i-\bar z\|^2,
\qquad
\bar z=\frac1p\sum_{i=1}^pz_i.
\end{equation}
Indeed, writing $\xi_i=\mathbf 1_{\{i\in J\}}$, we have
\[
\E[\xi_i]=\frac qp,
\qquad
\E[\xi_i\xi_j]
=
\frac{q(q-1)}{p(p-1)},
\qquad i\ne j,
\]
so \eqref{eq:finitePopulationVariance} follows by expanding the square and using $\sum_{i=1}^p(z_i-\bar z)=0$.
Applying \eqref{eq:finitePopulationVariance} for each $m$ and
using \eqref{eq:NystromCandidate} yields
\eqref{eq:NystromRandomRate}. 
Finally, since $\bar z_m$ is the mean of $z_{m,1},\ldots,z_{m,n_m}$,
\[
\frac1{n_m}\sum_{n=1}^{n_m}
\|z_{m,n}-\bar z_m\|_{\widehat T,\lambda}^2
=
\frac1{n_m}\sum_{n=1}^{n_m}
\|z_{m,n}\|_{\widehat T,\lambda}^2
-
\|\bar z_m\|_{\widehat T,\lambda}^2,
\]
which is exactly \eqref{eq:VmExplicit}.
\end{proof}

\begin{proof}[Proof of \cref{cor:NystromExplicit}]
By \eqref{eq:hLambda3}, we have $\widehat h_\lambda=\bm{A}^*\widehat{\bm\beta}_\lambda$ and
$(\bm{A}\bm{A}^*+\lambda\bm{G})
\widehat{\bm\beta}_\lambda
=
\bm y$,
whose $m$-th block reads
\(\lambda\frac{n_m}{w_m}\widehat{\beta}_m
=
\bm y_m-A_m\widehat h_\lambda,\) 
that is,
\[
n_m\widehat\beta_{m,n}
=
\frac{w_m}{\lambda}
\bigl(
f_m(x_{m,n})
-\calA_m\widehat h_\lambda(x_{m,n})
\bigr),
\qquad
n=1,\ldots,n_m.
\]
Dropping the nonpositive term in \eqref{eq:VmExplicit} therefore gives
\[
V_{\lambda,m}
\le
\frac1{n_m}
\sum_{n=1}^{n_m}
(n_m\widehat\beta_{m,n})^2
\|\varphi_{m,n}\|_{\widehat T,\lambda}^2
\le 
\frac{w_m}{\lambda^2}
\kappa_{\lambda,m,N}^2R_m^2.
\]
Inserting this bound into \eqref{eq:NystromRandomRate} yields the first inequality in \eqref{eq:NystromResidualRate}. The second follows by taking the maximum over $m$. 
\ref{cor:NystromExplicit1}: $R_N^2$ is the data-fit part of
$\widehat\Phi_\lambda(\widehat h_\lambda)$, and $\widehat h_\lambda$
minimizes $\widehat\Phi_\lambda$, whence
$R_N^2
\le
\widehat\Phi_\lambda(\widehat h_\lambda)
\le
\widehat\Phi_\lambda(0)$. If $u\in\calH$ solves \eqref{eq:generalProblem}, then $A_mu=\bm y_m$ for all $m$, and thus $\widehat\Phi_\lambda(\widehat h_\lambda)\le\widehat\Phi_\lambda(u)=\lambda\|u\|_\calH^2$. 
\ref{cor:NystromExplicit2}: By \eqref{eq:empiricalNorm} and the Cauchy--Schwarz inequality, with $(A_\ell\varphi_{m,n})_{n'}=k^{\ell,m}(x_{\ell,n'},x_{m,n})$ and $|k^{\ell,m}(x,y)|\le\kappa_\ell\kappa_m$,
\begin{align*}
\|\varphi_{m,n}\|_{\widehat T,\lambda}^2
=\sum_{\ell=1}^M\frac{w_\ell}{n_\ell}
\sum_{n'=1}^{n_\ell}
 k^{\ell,m}(x_{\ell,n'},x_{m,n})^2
 +\lambda\|\varphi_{m,n}\|_\calH^2\le\sum_{\ell=1}^M w_\ell\kappa_\ell^2\kappa_m^2
+\lambda\kappa_m^2.\qedhere
\end{align*}
\end{proof}

\begin{proof}[Proof of \cref{cor:NystromBalancedRate}]
The first inequality in \eqref{eq:NystromBalancedRate} is \eqref{eq:NystromResidualRate} with $n_m=n$ and $r_m=r_0$, and the second follows from $(n-r_0)/(n-1)\le1$ and $r_0=r/M$. The remaining statements follow from \cref{cor:NystromExplicit}\ref{cor:NystromExplicit1}--\ref{cor:NystromExplicit2} and  $\frac{w_m}{n_m}\|\bm y_m\|_2^2\le w_m\sup_{\calX_m}|f_m|^2$.
\end{proof}

\section{Additional Numerical Examples}\label{app:suppNumerics}

This appendix provides three additional numerical examples. We first consider a first-order ODE with the Gaussian kernel (\cref{ex:GaussianKernel}). We then introduce a finite-dimensional Fourier RKHS and use it to solve a second-order ODE. Finally, we illustrate a time-warped Gaussian kernel on the one-dimensional heat equation.

\subsection{First-Order ODE}\label{sec:supp1stODE}

Consider the first-order ODE
\begin{equation}\label{eq:supp1stODE}
    \begin{cases}
        u'(x)=-\pi\sin(\pi x), & x\in(-1,1),\\
        u(x)=-1, & x\in\{\pm1\},
    \end{cases}
\end{equation}
with $u'=\partial_1u$. Its unique solution is $u(x)=\cos(\pi x)$. While one of the two boundary conditions would suffice, we keep both to illustrate that the framework also handles overdetermined systems. We set $M=2$, $\calX_1=(-1,1)$, $\calX_2=\{\pm1\}$, and
\begin{align*}
    \calA_1&=\partial_1,
    & f_1(x)&=-\pi\sin(\pi x),\\
    \calA_2u&=u|_{\calX_2},
    & f_2&\equiv-1.
\end{align*}

\begin{figure}
    \centering
    \includegraphics[width=0.45\linewidth]{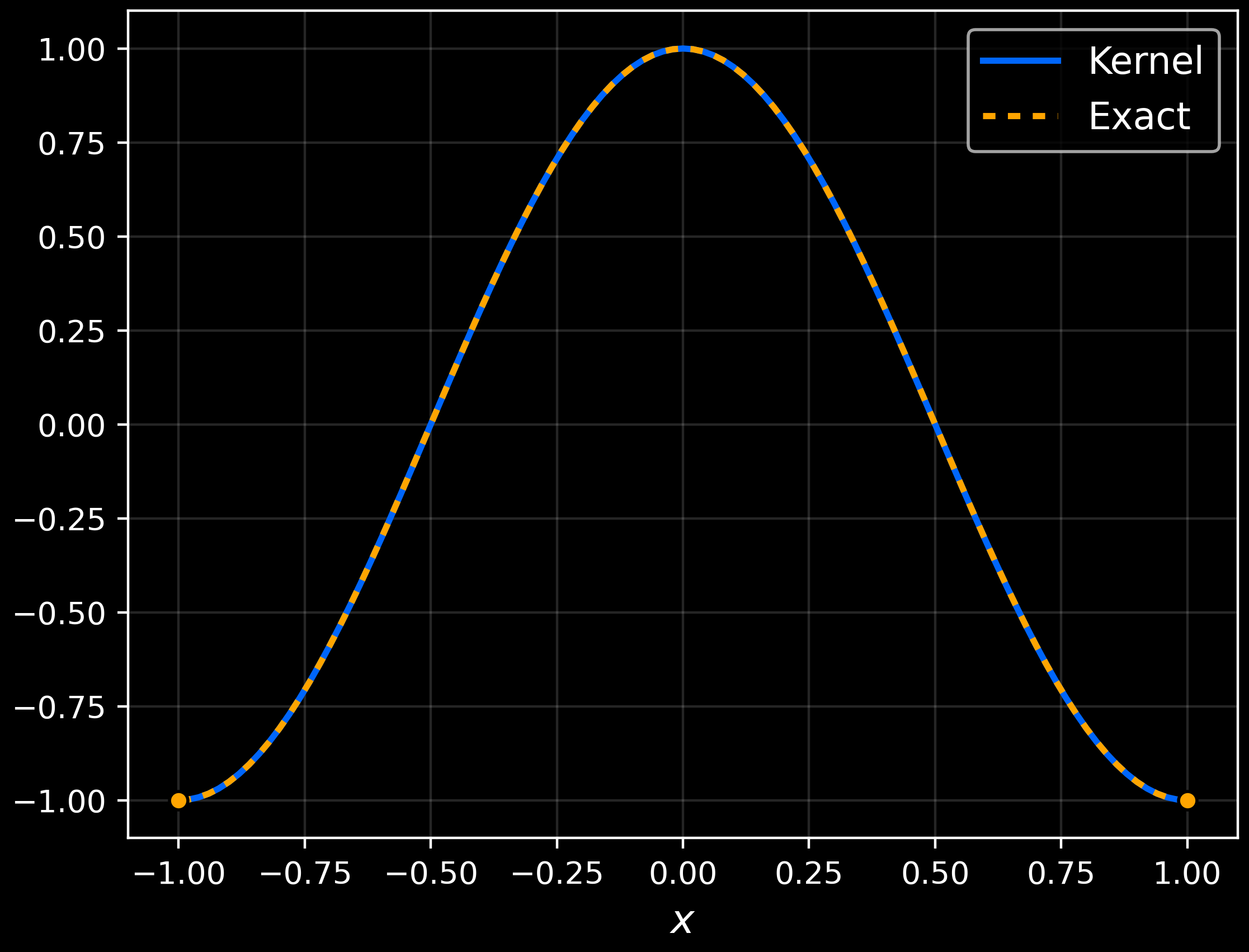}
    \caption{Kernel and exact solution of the first-order ODE \eqref{eq:supp1stODE}.}
    \label{fig:supp1stODE}
\end{figure}

We apply the kernel solver to \eqref{eq:supp1stODE} with regularization parameter $\lambda=10^{-8}$ and equal objective weights $w_1=w_2=1$. We use $n_1=100$ collocation points in $\calX_1$ drawn from the uniform distribution
\[
    \Pb_1(dx)=\frac12\mathbf 1_{(-1,1)}(x)\,dx.
\]
For the boundary condition, we directly use the $n_2=2$ deterministic points $x_{2,1}=-1$ and $x_{2,2}=1$, corresponding to the support of $\Pb_2=\frac12(\delta_{-1}+\delta_1)$. We use the Gaussian kernel $k(x,y)=e^{-\frac{\gamma}{2}|x-y|^2}$ with $\gamma=100$, corresponding to the length scale $\gamma^{-1/2}=0.1$. As shown in \cref{fig:supp1stODE}, the kernel solution $\widehat h_\lambda$ is visually indistinguishable from the exact solution.

\subsection{An ODE Example with a Finite-Dimensional Fourier RKHS}\label{sec:suppFourierODE}

We first introduce the Fourier kernel used in this example. Let $\calS=[-1,1]$, equipped with the Lebesgue measure. For a summable sequence $(\gamma_n)_{n\ge0}$ with $\gamma_n>0$, define
\begin{equation}\label{eq:suppFourierKernel}
    k(x,y)=\gamma_0+2\sum_{n=1}^\infty
    \gamma_n\cos\bigl(\pi n(x-y)\bigr).
\end{equation}
The functions
\[
\phi_0\equiv\frac{1}{\sqrt2},\qquad
\phi_n(x)=\cos(\pi n x),\qquad
\psi_n(x)=\sin(\pi n x),\quad n\ge1,
\]
form an orthonormal basis of $L^2(\calS)$. Moreover,
\[
k(x,y)
=2\gamma_0\phi_0(x)\phi_0(y)
+2\sum_{n=1}^\infty\gamma_n
\bigl(\phi_n(x)\phi_n(y)+\psi_n(x)\psi_n(y)\bigr).
\]
Hence the associated RKHS consists of all functions
\begin{align*}
h&=h_0\phi_0+\sum_{n\ge1}(h_n^\phi\phi_n+h_n^\psi\psi_n),\\
h_0&=(h,\phi_0)_{L^2},\qquad
h_n^\phi=(h,\phi_n)_{L^2},\qquad
h_n^\psi=(h,\psi_n)_{L^2},
\end{align*}
satisfying
\[
\|h\|_{\calH}^2
=\frac{h_0^2}{2\gamma_0}
+\frac12\sum_{n=1}^\infty
\frac{(h_n^\phi)^2+(h_n^\psi)^2}{\gamma_n}
<\infty.
\]
If $\sum_{n\ge1} n^2\gamma_n<\infty$, the series can be differentiated termwise, and the associated RKHS $\calH$ is $\partial_1$-regular on $\calX=(-1,1)$, with
\[
\|k_x^{\partial_1}\|_\calH^2
=
2\sum_{n\ge1}\gamma_n(\pi n)^2.
\]
Similarly, $\calH$ is $\partial_{11}$-regular if $\sum_{n\ge1}n^4\gamma_n<\infty$, with
\[
\|k_x^{\partial_{11}}\|_\calH^2
=2\sum_{n\ge1}\gamma_n(\pi n)^4.
\]
The decay of the weights $\gamma_n$ thus determines for which differential operators $\calH$ is regular.

For numerical purposes, we replace the infinite series by its finite truncation at level $L\in\N$, taking $\gamma_0=1$, $\gamma_n=1/(2n^4)$ for $1\le n\le L$, and $\gamma_n=0$ otherwise. This yields
\begin{align}
    k_L(x,y)
    &=1+\sum_{n=1}^L\frac{1}{n^4}\cos\bigl(\pi n(x-y)\bigr),
    \label{eq:suppFourierTrunc}
    \calH_L
    &=\operatorname{span}\{\phi_0,\phi_n,\psi_n:1\le n\le L\}.
\end{align}
The inner product is
\[
(h,g)_{\calH_L}
=\frac12h_0g_0+\sum_{n=1}^L
n^4\bigl(h_n^\phi g_n^\phi+h_n^\psi g_n^\psi\bigr).
\]
Hence
\[
\|h\|_{\calH_L}^2
=\frac12h_0^2+\pi^{-4}\|h''\|_{L^2(\calS)}^2,
\]
so that the weights penalize high frequencies as a Sobolev norm of order two. The truncated kernel is smooth, and therefore $\calH_L$ is $\partial_1$- and $\partial_{11}$-regular.

\begin{figure}
    \centering
    \includegraphics[width=0.95\linewidth]{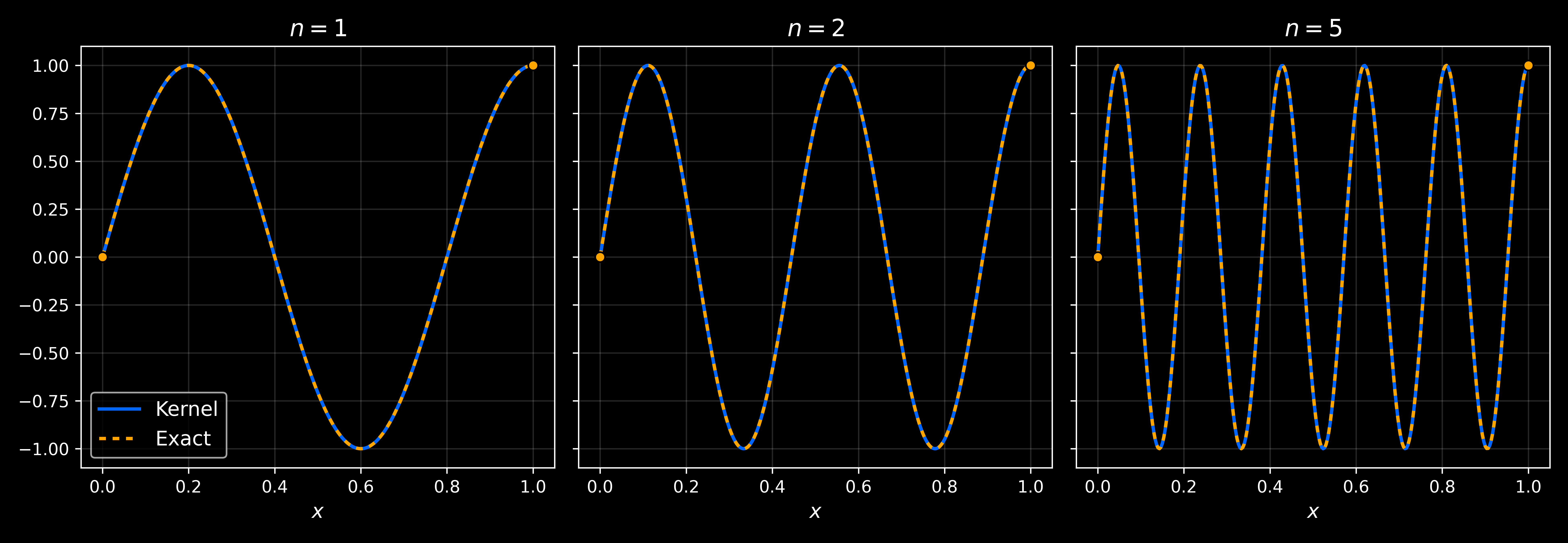}
    \caption{Kernel and exact solution of the second-order ODE \eqref{eq:supp2ndODE} for $n\in\{1,2,5\}$.}
    \label{fig:supp2ndODE}
\end{figure}

We now use this kernel to solve a second-order ODE. Let
$\zeta_n=(2n+1/2)^2\pi^2$, $n\in\N$, and consider
\begin{equation}\label{eq:supp2ndODE}
    \begin{cases}
        u''+\zeta_nu=0, & \textnormal{on }\calX_1=(0,1),\\
        u(x)=x, & \textnormal{on }\calX_2=\{0,1\},
    \end{cases}
\end{equation}
where $u''=\partial_{11}u$. The unique solution of \eqref{eq:supp2ndODE} is
\[
    u_n(x)=\sin(\sqrt{\zeta_n}\,x).
\]
We use the truncated Fourier kernel \eqref{eq:suppFourierTrunc}, evaluated on $[0,1]$, with truncation level $L=20$. We take $\lambda=10^{-10}$, equal objective weights $w_1=w_2=1$, $n_1=150$ equally spaced collocation points in $[0.01,0.99]\subset\calX_1$, and the $n_2=2$ boundary points $x\in\{0,1\}$. The results are shown in \cref{fig:supp2ndODE} for $n\in\{1,2,5\}$, where the kernel and exact solutions closely agree.

\begin{remark}
Observe that $u_n$ is the $(2n+1)$-th eigenfunction of the covariance kernel of Brownian motion on $[0,1]$. That is, $u_n$ solves the following Fredholm integral equation of the second kind 
\[
    \int_0^1 u_n(y)g(x,y)\,dy
    =\zeta_n^{-1}u_n(x),
    \qquad
    g(x,y)=x\wedge y.
\]
The integral equation alone determines the eigenfunction up to a multiplicative constant. In particular, it admits the zero solution. Together with the normalization $u_n(1)=1$ from \eqref{eq:supp2ndODE}, it determines $u_n$ uniquely. The kernel solver applies to this equivalent Fredholm formulation as well, using the integral-operator construction of \cref{ex:integralEquations} together with the stability of operator representers under linear combinations in \cref{eq:lincomb}. We do not pursue this here. 
\end{remark}

\subsection{Heat Equation}\label{sec:suppHeat}

For this forward parabolic example, we use a Gaussian kernel whose spatial coordinate is warped according to the diffusion scaling. Given $(t,x),(s,y)\in\calS\coloneqq[0,\infty)\times\R$, define
\begin{align}
k\big((t,x),(s,y)\big)
&=
\exp\left\{
-\frac{\gamma_0}{2}|t-s|^2
-\frac{\gamma_1}{2}|\psi(t,x)-\psi(s,y)|^2
\right\},
\label{eq:suppFwdParabolicKernel}\\
\psi(t,x)&=\frac{x}{\sqrt{t+\varepsilon}},
\qquad \varepsilon>0.
\label{eq:suppFwdTransform}
\end{align}
The transformation reflects the parabolic scaling $x\sim\sqrt{t}$. In particular, at fixed $t$, the spatial length scale grows proportionally to $\sqrt{t+\varepsilon}$. The parameter $\varepsilon>0$ prevents a singularity at $t=0$ and ensures that the kernel remains smooth.

\begin{figure}[H]
    \centering
    \begin{subfigure}[t]{0.48\linewidth}
        \centering
        \includegraphics[width=0.98\linewidth]{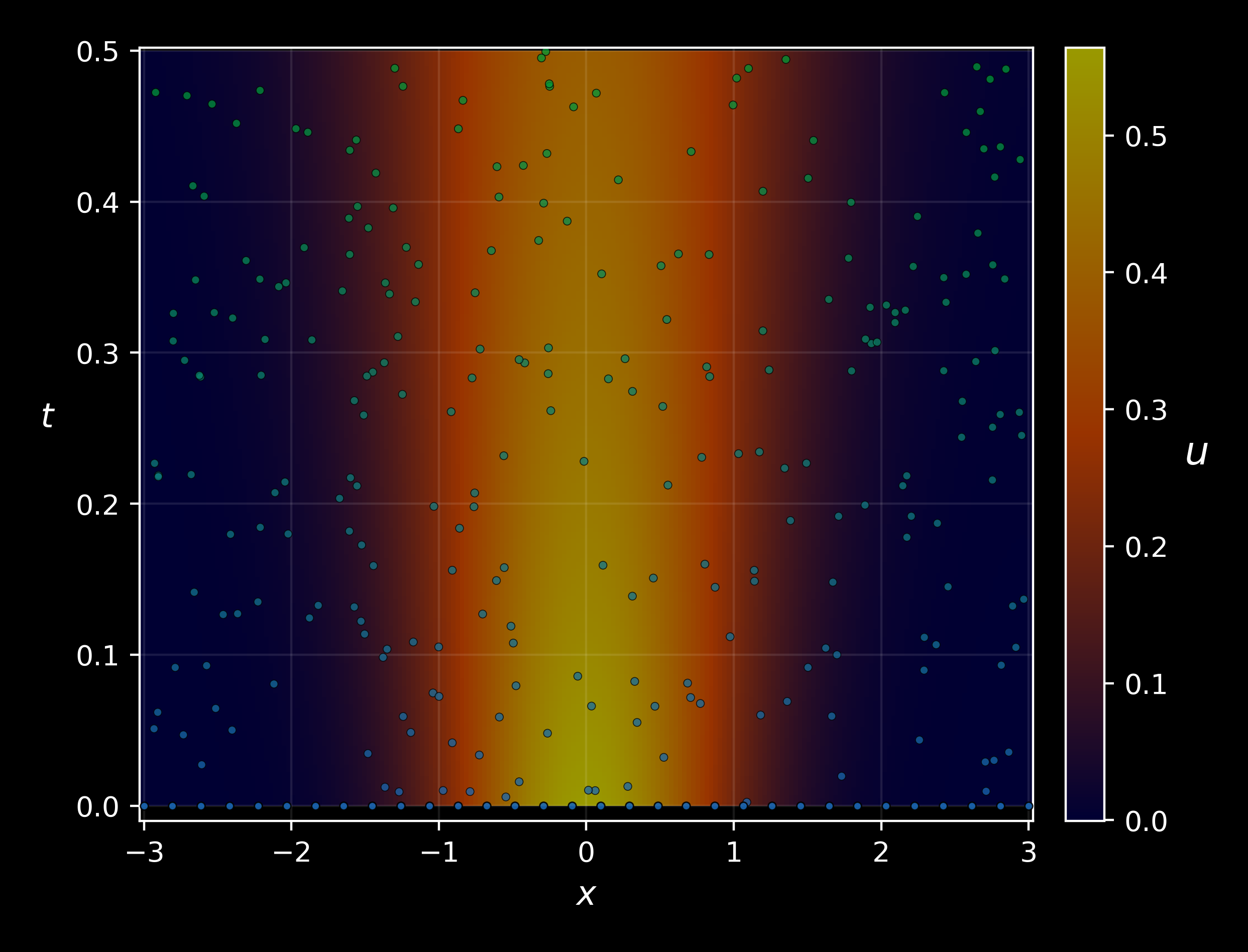}
        \label{fig:suppHeatSamples}
    \end{subfigure}
    \hfill
    \begin{subfigure}[t]{0.48\linewidth}
        \centering
        \includegraphics[width=0.98\linewidth]{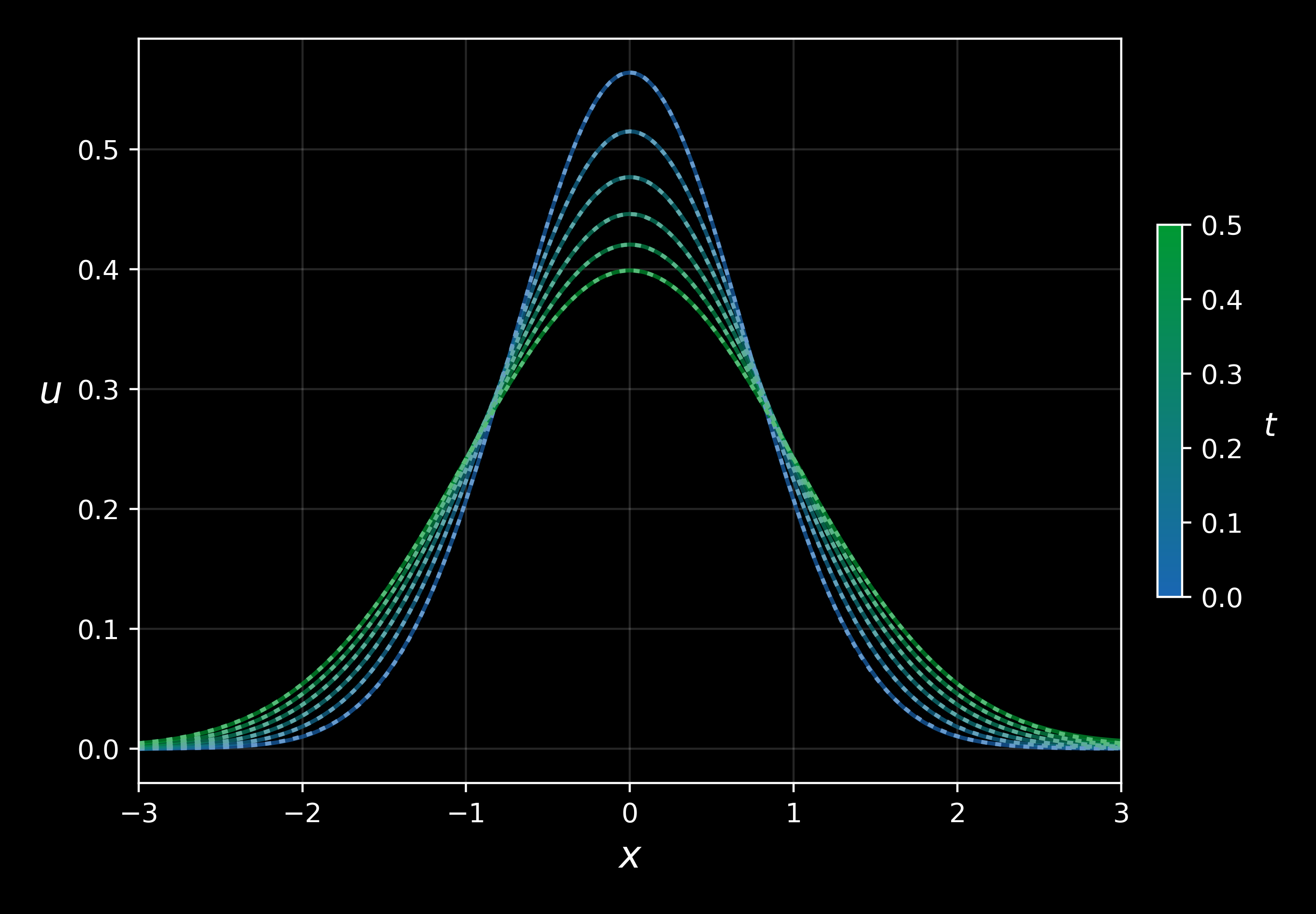}
        \label{fig:suppHeatSolution}
    \end{subfigure}
    \caption{One-dimensional heat equation \eqref{eq:supp1DHeat}. Left: parabolic domain, collocation points, and kernel solution. Right: time sections of the kernel solution (solid lines) compared with the exact solution (dashed lines).}
    \label{fig:suppHeatEquation}
\end{figure}

Consider the one-dimensional heat equation
\begin{equation}\label{eq:supp1DHeat}
    \begin{cases}
        \bigl(\partial_t-\frac12\partial_{xx}\bigr)u=0,
        & \textnormal{on }\calX_1=(0,T)\times\R,\\
        u(0,x)=\phi_{1/2}(x),
        & \textnormal{on }\calX_2=\{0\}\times\R,
    \end{cases}
\end{equation}
where $\phi_t$ denotes the density of $\calN(0,t)$. We take $T=0.5$ and the spatial box $[-3,3]$. We use $n_1=256$ collocation points drawn uniformly in $(0,T)\times[-3,3]$ and $n_2=32$ equally spaced collocation points in $\{0\}\times[-3,3]$, see the left panel of \cref{fig:suppHeatEquation}. We use the parabolic kernel \eqref{eq:suppFwdParabolicKernel} with $\gamma_0=4$, $\gamma_1=1$, and $\varepsilon=0.25$. Thus the time bandwidth is $\gamma_0^{-1/2}=0.5$ and the bandwidth in the transformed spatial coordinate is $\gamma_1^{-1/2}=1$. We use equal objective weights $w_1=w_2=1$ and Tikhonov regularization parameter $\lambda=10^{-7}$.

The right panel of \cref{fig:suppHeatEquation} compares the kernel solution $\widehat h_\lambda$ with the exact solution $u(t,x)=\phi_{1/2+t}(x)$. The two functions are visually indistinguishable. The relative $L^2$ error on a regular grid of $[0,T]\times[-3,3]$ is about $1.3\cdot10^{-3}$.

\bibliographystyle{abbrvnat}
\bibliography{ref}

\end{document}